\documentclass{amsart}

\usepackage{amssymb,amsthm}
\usepackage{eucal}

\def\cl {L}

\def\ce {\mathcal{E}}
\def\ca {\mathcal{A}}
\def\ci {\mathcal{I}}
\def\cb {\mathcal{B}}

\def\fm {\mathfrak{m}}
\def\fM {\mathfrak{M}}
\def\fs {\mathfrak{S}}
\def\fz {\mathfrak{Z}}
\def\beq {\begin{equation}}
\def\endq {\end{equation}}
\newcommand{\twolinesum}[2]{\sum_{\substack{{\scriptstyle #1}\\
{\scriptstyle #2}}}}
\newcommand{\threelinesum}[3]{\sum_{\substack{{\scriptstyle #1}\\
{\scriptstyle #2}\\{\scriptstyle #3}}}}

\newcommand{\modulo}[1]{\mathrm{mod}\;#1}
\newcommand{\pmodulo}[1]{\;(\mathrm{mod}\;#1)}
\newcommand{\sfrac}[2]{{\textstyle \frac{#1}{#2}}}

\newtheorem{theorem}{Theorem}
\newtheorem{lemma}{Lemma}

\newtheorem*{theorem*}{Theorem}

\theoremstyle{remark}
\newtheorem*{remark}{Remark}

\numberwithin{equation}{section}

\begin{document}

\title{On sums of squares of primes II}
\author{Glyn Harman}
\address{Department of Mathematics\\
         Royal Holloway University of London\\
         Egham\\
         Surrey TW20 0EX\\
         U.K.}
\email{G.Harman@rhul.ac.uk}

\author{Angel Kumchev}
\address{Department of Mathematics\\
         Towson University\\
         7800 York Road\\
         Towson, MD 21252\\
         U.S.A.}
\email{akumchev@towson.edu}


\begin{abstract}
In this paper we continue our study, begun in \cite{hk1}, of the exceptional set of integers,
not restricted by elementary congruence conditions, which cannot
be represented as sums of three or four squares of primes.  We correct a serious oversight in our first paper,
but make further progress on the exponential sums estimates needed, together with an embellishment of
the previous sieve technique employed.  This leads to an improvement in our bounds for the
maximal size of the exceptional sets.
\end{abstract}

\maketitle

\section{Introduction}

As in \cite{hk1} we write:
\[
\begin{split}
\ca_3 &= \{n \in \mathbb{N}: n \equiv 3 \pmodulo {24}, \;
           n \not\equiv 0 \pmodulo 5 \},\\
\ca_4 &= \{n \in \mathbb{N}: n \equiv 4 \pmodulo {24} \}.
\end{split}
\]
We further put:
\[
E_j(N) = |\{n \in \ca_j: n \le N, n \ne p_1^2 + \cdots + p_j^2, \text{for any primes $p_u$}  \}|, \ \ j=3,4.
\]
Our purpose in writing this article is to correct an error in our previous discussion of upper bounds for these sets
and also to introduce further refinements to the method which lead to superior results.  Although the improvement
in the exponent is relatively small (the crucial change is from $1/7$ to $3/20$) the modifications to the method have independent interest
and may have further applications - we state one such result below as Theorem 3.

It is conjectured that every sufficiently large integer in
$\ca_j$ can be represented as the sum of $j$ squares of primes, and so $E_j(N) = O(1)$.  The expected main terms from
an application of the Circle Method lead one to the following
hypothetical asymptotic formulae:
\beq\label{f1}
  \sum_{p_1^2 + p_2^2 + p_3^2 = n}
  (\log p_1)(\log p_2)(\log p_3)
  \ \sim \ \frac{\pi}{4} \fs_3(n) n^{1/2}
\endq
and
\beq\label{f2}
  \sum_{p_1^2 + \dots + p_4^2 = n}
  (\log p_1) \cdots (\log p_4)
  \ \sim \ \frac{\pi^2}{16} \fs_4(n) n,
\endq
where $\fs_j(n) > 0$ for all large $n \in \ca_j$.  In 1938  Hua \cite{Hua} proved a general result on representing almost all numbers in
suitable residue classes as the sum of two squares of primes and the $k-$th power of a prime, from which it follows
 that almost all $n \in \mathcal A_3$ are representable as sums of
three squares of primes. Of course, we then immediately obtain that almost all
$n \in \mathcal A_4$ are representable as sums of four squares of primes.
The subsequent history of this problem is documented in \cite{hk1} (charting the developments in \cite{Schwarz,Leung,Zhan,LiuWooley,Kumchev}), culminating in
the authors' demonstration that
\[
  E_3(N) \ll N^{{6}/{7} + \epsilon} \qquad \text{and} \qquad
  E_4(N) \ll N^{5/14 + \epsilon}.
\]
Unfortunately there was a serious oversight in our proofs. To be precise, the display (4.16) in \cite{hk1}
which gives an estimate on average for the singular series, namely
\[
\sum_{N/2 < n \le N} \bigg| \fs_3(n, Q) -
  8\prod_{2 < p \le Q} (1 + \mathfrak s(p, n)) \bigg|
  \ll N^{1 + \epsilon/2}Q^{-1/2},
\]
is not good enough for the stated result to follow.
We would like to thank Claus Bauer and Hongze Li who both independently alerted the authors to this error.
We correct this in section 3.5 here.  However, we can now move beyond what seemed a difficult barrier
with the previous exponents, which arose as $1 - \sigma$ and $\frac12 - \sigma$ with $\sigma = 1/7$.
The following results show that we can now increase $\sigma$ to $3/20$.

\begin{theorem}\label{th1}
Let $\epsilon > 0$ be given.  Then for all large $N$ we have
\beq\label{f3}
E_3(N) \ll N^{17/20 + \epsilon}.
\endq
\end{theorem}

\begin{theorem}\label{th2}
Let $\epsilon > 0$ be given.  Then for all large $N$ we have
\beq\label{f4}
E_4(N) \ll N^{7/20 + \epsilon}.
\endq
\end{theorem}

Combining the new ideas in the present work with \cite{HL} we obtain the following.

\begin{theorem}\label{th3}
Let $E(N)$ represent the cardinality of the set
\[
\{n \le N: n \equiv 1 \ \text{or} \ 3 \pmodulo{6}, n \ne p_1 + p_2^2 + p_3^2   \}.
\]
Then, for every $\epsilon > 0$,
\beq\label{HLimproved}
E(N) \ll N^{7/20 + \epsilon}.
\endq
\end{theorem}

\section{The Method}

We shall only prove Theorem \ref{th1}; the straightforward modifications needed for Theorem
\ref{th2} follow as in \cite{hk1}, and for Theorem \ref{th3} as in \cite{HL}. It suffices to estimate the
number of exceptional integers $n$ in the set $\cb = \ca_3
\cap (\frac12 N, N]$ where $N$ will be our main parameter, which we
assume to be ``sufficiently large". We write
\[
  P = N^{1/2}, \qquad \cl = \log P, \qquad
  \ci = \left[ \textstyle \frac 13 P, \frac 23 P \right).
\]
We use $c$ to denote an absolute constant, not necessarily the
same at each occurrence. In the following, $\sigma$ will be a parameter in the range
$\frac17 \le \sigma \le \frac{3}{20}$, and our method will show that
$E_3(N) \ll N^{1 - \sigma + \epsilon}$. Here, as elsewhere in the following, $\epsilon$ is
an arbitrary small positive real.

We wish to represent integers $n$ in the form $m_1^2 +m_2^2 +m_3^2$
where each $m_j$ is restricted to prime values.
In our previous paper we sieved only one of the variables, say $m_3$.
In our current work we will sieve two variables, albeit in a rather asymmetric way.
To be precise, let $\rho_1(m)$ be the characteristic function of the set of primes.
Suppose that, for suitable non-negative functions $\rho_j(m), 2 \le j \le 5$, we have
\[
\rho_2(m) \ge \rho_1(m) = \rho_3(m) - \rho_4(m) + \rho_5(m).
\]
Then
\[
\twolinesum{m_1^2 +m_2^2 +m_3^2 = n}{m_j \in \ci} \rho_1(m_1)\rho_1(m_2)\rho_1(m_3)
\ge
S_1 - S_2
\]
where
\[
\begin{split}
S_1 &= \twolinesum{m_1^2 +m_2^2 +m_3^2 = n}{m_j \in \ci} \rho_1(m_1)\rho_1(m_2)\rho_3(m_3), \\
S_2 &= \twolinesum{m_1^2 +m_2^2 +m_3^2 = n}{m_j \in \ci} \rho_1(m_1)\rho_2(m_2)\rho_4(m_3).
\end{split}
\]

The circle
method then gives
\beq\label{HL}
 \twolinesum{m_1^2 +m_2^2 +m_3^2 = n}{m_j \in \ci} \rho_j(m_1)\rho_k(m_2)\rho_{\ell}(m_3)
  =
  \int_0^1 f_j(\alpha) f_k(\alpha) f_{\ell}(\alpha) e(-\alpha n) \, d\alpha,
\endq
where we write $e(x) = \exp(2 \pi \mathrm i x)$ and, for $1 \le j \le 4$,
\beq\label{expsums}
f_j(\alpha) = \sum_{m \in \ci} \rho_j(m) e(\alpha m^2).
\endq

Here we will want the $\rho_j, 2 \le j \le 4$, to satisfy:
\beq\label{rho}
\sum_{m \le X} \rho_j(m) =C_j  X\cl^{-1}(1 + o(1))
\endq
for $P^{1/2} \le X \le P$, where
\beq\label{constants}
C_3 - C_2C_4 > 0.
\endq
It then remains to establish that
\[
\int_0^1 f_1(\alpha) f_k(\alpha) f_{\ell}(\alpha) e(-\alpha n) \, d\alpha = K_n C_k C_{\ell}
 \Pi(n,Q) P L^{-3}(1 + o(1))
\]
for the same value $K_n$ in the two cases $k=1, \ell=3$, $k=2, \ell=4$ where $C_1 = 1$,
with at most $E_3(N)$ exceptions up to $N$. Here $\Pi(n,Q)$ is an approximation
to $\fs(n)$ which we define later and which satisfies $\Pi(n,Q) \gg L^{-3}$. When we state the main term more
explicitly it will be clear that $1 \ll K_n \ll 1$ with absolute constants.
The properties of the $\rho_j$ necessary to achieve this will be
introduced when relevant.  In particular it should be noted that we require
$\rho_2$ and $\rho_3$ to satisfy the most stringent conditions.

Our application of the circle method has the same format as our previous work; see
 \cite{Vaughan} for a general introduction. The main contribution to the right side of
\eqref{HL} comes from the \emph{major arcs} which we denote by $\fM$ and are defined as follows.  Let
$Q = P^{2 \sigma - 3\epsilon}$ and write (shifting
$[0,1)$ by $\omega = QP^{-2 + \epsilon}$ which does not change \eqref{HL})
\beq\label{arcs1}
\fM = \left[\omega, 1 + \omega \right)
\cap \bigcup_{1 \le q \le Q} \bigcup_{(a,q)= 1} \left[\frac{a}{q} - \frac{\omega}{q},
\frac{a}{q} + \frac{\omega}{q} \right).
\endq
The \emph{minor arcs} $\fm$ are then given by
$\fm = [\omega, 1 + \omega) \setminus \fM.$

For technical reasons, it is convenient to modify $f_j(\alpha)$, $j \ge 2$,  on the major arcs to remove
interference between possible prime divisors of $m$ (when $\rho(m) < 0$) and approximation denominators.
We introduce a function $\theta(m, \alpha)$ which is $1$ except when there exist integers $a$ and $q$ such that
\[
  |q\alpha - a| < \omega, \quad (a,q) = 1, \quad  q \le Q, \quad (m, q) \ge P^{\sigma},
\]
in which case $\theta(m, \alpha) = 0$. Write
\[
  g_j(\alpha) = \sum_{m \in \ci} \rho_j(m) \theta(m, \alpha)e(\alpha m^2).
\]
We note that $g_j(\alpha) = f_j(\alpha)$ for $\alpha \in \fm$ and that
\beq\label{2.5}
  f_j(\alpha) - g_j(\alpha) \ll P^{1 - \sigma}
\endq
for all $\alpha$.

\section{The major arcs}
\label{secM}

The major arc contributions to $S_1$ and $S_2$ are dominated by the integrals
\[
  \int_{\fM} f_1(\alpha)^2g_3(\alpha)e(-\alpha n) \, d\alpha \quad \text{and} \quad
  \int_{\fM} f_1(\alpha)g_2(\alpha)g_4(\alpha)e(-\alpha n) \, d\alpha,
\]
respectively. In this section, we evaluate the latter integral. The evaluation of the former can be carried out
in a similar fashion and is, in fact, less technical.

As in \cite{hk1}, we suppose that $\rho_j$, $j = 2, 3, 4$, have asymptotic properties similar to those of $\rho_1$.
To be precise, we assume that $\rho_j$ satisfy the following two hypotheses:
\begin{itemize}
  \item [(i)] Let $A, B > 0$ be fixed, let $\chi$ be a non-principal character modulo $q$, $q \le \cl^B$, and let $\ci'$ be a subinterval of $\ci$. Then
  \begin{equation}\label{M.1}
    \sum_{m \in \ci'} \rho_j(m) \chi(m) \ll P\cl^{-A}.
  \end{equation}
  \item [(ii)] Let $A > 0$ be fixed and let $\ci'$ be a subinterval of $\mathcal{I}$. There exists a smooth function $\delta_j$ on $\ci$ such that
  \begin{equation}\label{M.2}
    \sum_{m \in \ci'} \rho_j(m) = \sum_{m \in \ci'} \delta_j(m) + O \big( P\cl^{-A} \big).
  \end{equation}
\end{itemize}
Of course, by the Siegel--Walfisz theorem, these hypotheses hold also for $\rho_1(m)$ with $\delta_1(m) = (\log m)^{-1}$.
We note that \eqref{M.2} gives
\[
\int_{\ci} \delta_j(u) \, du = C_j \frac{P}{3L}(1 + o(1)).
\]
 Furthermore, we assume that:
\begin{itemize}
  \item [(iii)] $\rho_j(m) = 0$ if $m$ has a prime divisor $p < Z = P^{1 - 6\sigma}$.
\end{itemize}

For $j = 1, \dots, 4$, we define functions $f_j^*(\alpha)$ on $\fM$ by setting
\[
  f_j^*(\alpha) = \frac {S(\chi_0, a)}{\phi(q)} \sum_{m \in \ci} \delta_j(m) e \left( (\alpha - a/q) m^2 \right) \quad \text{if } \alpha \in \fM(q, a).
\]
Here $\chi_0$ is the principal character modulo $q$ and
\[
  S(\chi, a) = \sum_{h = 1}^q \bar{\chi}(h)e_q(ah^2).
\]
We now proceed to estimate the integral
\begin{equation}\label{M.3}
  \int_{\fM} \big( f_1(\alpha)g_2(\alpha)g_4(\alpha) - f_1^*(\alpha)f_2^*(\alpha)f_4^*(\alpha) \big) e(-\alpha n) \, d\alpha,
\end{equation}
which we think of as the error of approximation of the contribution from $\fM$ by the expected main term.
For our purposes, it suffices to show that this quantity is $O(P\cl^{-A})$ for any fixed $A > 0$, for example.

A difficulty arises upon reducing $\sigma$ below $1/7$ -- the function $\theta(m,\alpha)$ no longer covers
the interference between all possible prime divisors of $m$ (when $\rho(m) < 0$) and the major arc denominators. To be precise, we need a new argument
for the range from $Z$ to $P^{\sigma}$.
To deal with this,
for an integer $q$, we write $\mathcal S_q$ for the set of primes $p$ in the range $Z \le p < P^{\sigma}$ that divide $q$.
In particular, $\mathcal S_0$ is simply the set of primes $p$ with $Z \le p < P^{\sigma}$. We also write $\mathcal S_q' = \mathcal S_q \cup \{1\}$.
Since $Z^2 > P^{\sigma}$, under hypothesis (iii), we have
\[
  g_j(\alpha)  = \sum_{l \in \mathcal S_q'} g_{j,l}(\alpha) = g_{j,1}(\alpha) + \sum_{p \in \mathcal S_q} g_{j,p}(\alpha),
\]
where for $\alpha \in \fM(q, a)$ and $l \in \mathcal S_q'$,
\[
  g_{j,l}(\alpha) = \twolinesum{m \in \ci}{(m,q) = l} \rho_j(m) \theta(m, \alpha) e(\alpha m^2).
\]
Similarly to (4.1) in \cite{hk1}, when $\alpha \in \mathfrak M(q, a)$ and $l \in \mathcal S_q'$, we have
\begin{equation}\label{M.4}
  g_{j,l}(\alpha) = \frac 1{\phi(q_l)} \sum_{\chi \bmod q_l} S(\chi, al)\sum_{lm \in \ci} \rho_j(lm)\chi(m) e(\beta l^2m^2),
\end{equation}
where $q_l = q/l$ and $\beta = \alpha - a/q$. If $\chi$ is a character and $l$ a natural number, we now define
\[
  W_{j,l}(\chi, \beta) = \sum_{lm \in \ci} (\rho_j(lm) \chi(m) - D_l(\chi)\delta_j(lm)) e(\beta l^2m^2),
\]
where $D_l(\chi) = 1$ when $l = 1$ and $\chi$ is principal and $D_l(\chi) = 0$ otherwise. By \eqref{M.4} above and (4.1) in \cite{hk1},
\begin{align}
  \Delta_1(\alpha) &= f_1(\alpha) - f_1^*(\alpha) = \frac 1{\phi(q)} \sum_{\chi \bmod q} S(\chi, a)W_{1,1}(\chi, \alpha - a/q), \label{M.5}\\
  \Delta_j(\alpha) &= g_{j,1}(\alpha) - f_j^*(\alpha) = \frac 1{\phi(q)} \sum_{\chi \bmod q} S(\chi, a)W_{j,1}(\chi, \alpha - a/q), \label{M.6}\\
  g_{j,p}(\alpha) &= \frac 1{\phi(q_p)} \sum_{\chi \bmod q_p} S(\chi, ap)W_{j,p}(\chi, \alpha - a/q). \label{M.7}
\end{align}
Using \eqref{M.5}--\eqref{M.7}, we can express the integral \eqref{M.3} as the linear combination of seventeen quantities of the form
\[
  \int_{\mathfrak M} \Delta_1^{\flat}(\alpha)\Delta_2^{\flat}(\alpha)\Delta_4^{\flat}(\alpha)e(-\alpha n) \, d\alpha,
\]
with $\Delta_1^{\flat}(\alpha)$ one of $f_1^*(\alpha)$ or $\Delta_1(\alpha)$ and $\Delta_j^{\flat}(\alpha)$, $j = 2, 4$, one of
\[
  f_j^*(\alpha), \quad \Delta_j(\alpha) \quad \text{or} \quad \sum_{p \in \mathcal S_q} g_{j, p}(\alpha).
\]
To be more precise, each of the eighteen possible combinations occurs with the exception of
$f_1^*(\alpha) f_2^*(\alpha) f_3^*(\alpha)$ which we later show to give the
main term.

We shall restrict our attention here to the two most troublesome combinations:
\begin{gather}
  I_1 = \int_{\mathfrak M} \Delta_1(\alpha)\Delta_2(\alpha)\Delta_4(\alpha)e(-\alpha n) \, d\alpha, \label{M.8} \\
  I_2 = \sum_{p_1, p_2 \in \mathcal S_0} \int_{\mathfrak M_{\mathbf p}} \Delta_1(\alpha)g_{2, p_1}(\alpha)g_{4, p_2}(\alpha)e(-\alpha n) \, d\alpha, \label{M.9}
\end{gather}
where $\mathfrak M_{\mathbf p}$ denotes the subset of $\mathfrak M$ consisting of the major arcs $\mathfrak M(q, a)$, with $q$ divisible by $p_1$ and $p_2$.
However, before we estimate $I_1$ and $I_2$, we need to establish some lemmas.

\subsection{Bounds for averages of $W_{j,l}(\chi, \beta)$}
\label{secM.1}

At this point, we need to make a hypothesis about the structure of the sieve weights $\rho_j$. Henceforth, we write
\begin{equation}\label{M.10}
  \psi(m, z) = \begin{cases}
    1  & \text{if }  p \mid m \Rightarrow p \ge z, \\
    0  & \text{otherwise.}
  \end{cases}
\end{equation}
We also extend $\psi(m, z)$ to all real $m > 0$ by setting $\psi(m, z) = 0$ when $m$ is not an integer.
Our construction will yield coefficients $\rho_j$ that are linear combinations of convolutions of the form
\begin{equation}\label{M.11}
  \sum_{r \sim R} \sum_{s \sim S} \xi_r\eta_s\psi(rs, z)\psi(m/rs, z),
\end{equation}
where $|\xi_r| \le \tau(r)^c$ and $|\eta_s| \le \tau(s)^c$. In our applications the value of $z$ will often depend on certain variables.
To help set up the necessary hypotheses for our auxiliary results we therefore write $z(r,s)$ for a positive real-valued function,
which in practice will either be fixed, or take the value $p$ for some prime divisor of $r$ or $s$; see \S \ref{sec5} for the
specific cases of interest. We also put
\begin{equation}\label{M.12}
  Y = P^{1 - 5\sigma}, \quad V = P^{2\sigma}, \quad W = P^{1 - 4\sigma}.
\end{equation}
We now require that $\rho_j$ satisfies the following additional hypothesis:
\begin{itemize}
  \item [(iv)] $\rho_j$ can be expressed as a linear combination of $O(\cl^c)$ convolutions of the form \eqref{M.11}, where
    \begin{equation}\label{M.15}
    1 \le R \le V, \quad 1 \le S \le W, \quad Z \le z(r, s) \le P^{8/35}.
  \end{equation}
\end{itemize}

For the remainder of \S \ref{secM.1}, we suppress the index $j$ and write $W_l(\chi, \beta)$ for $W_{j,l}(\chi, \beta)$, $\rho$ for $\rho_j$, etc.

\begin{lemma}\label{lM.0}
  Let $\alpha, \beta$ be reals with $0 < \alpha < \beta$, let $n, g$ be positive integers, and let $(A_q)$ be a sequence of positive reals such that
  \[
    \twolinesum{q \sim Q}{d \mid q} A_q \le B_1 + d^{-1}B_2.
  \]
  Then
  \[
    \sum_{q \sim Q} (n, [q, g])^{\alpha}[q, g]^{-\beta}A_q \ll (n, g)^{\alpha}g^{-\beta + \epsilon}\big( B_1 + Q^{-\beta'}B_2 \big),
  \]
  where $\beta' = \min(\beta - \alpha, 1)$. Furthermore, if $ghQ \ge n^{\delta}$ for some $\delta > 0$, then
  \[
    \sum_{q \sim Q} (n, [q, g])^{\alpha}[q, g]^{-\beta}A_q \ll (n, g)^{\alpha}g^{-\beta + \epsilon}\big( B_1 + Q^{-\beta''}B_2 \big),
  \]
  where $\beta'' = \min(\beta, 1)$.
\end{lemma}

\begin{proof}
  These inequalities can be established by a slight generalization of the arguments leading to (5.21) and (5.23) in \cite{LiTs05}.
In particular, see (5.20) and (5.22) in \cite{LiTs05}.
\end{proof}

\begin{lemma}\label{lM.1}
  Suppose that $\rho$ is a convolution of the form \eqref{M.11} and $\Phi$ is a complex-valued function defined on $\ci$.
Suppose also that the parameters $R$ and $S$ and the function $z(r, s)$ satisfy
  \begin{equation}\label{M.13}
    \max(R, S) \le P^{11/20}, \quad z(r, s)\min(R, S) \le P^{11/20}, \quad z(r, s) \le P^{8/35}.
  \end{equation}
  Then the sum
  \[
    \sum_{m \in \ci} \rho(m)\Phi(m)
  \]
  can be expressed as a linear combination of $O(\cl^c)$ sums of the form
  \begin{equation}\label{M.14}
    \sum_{r \sim R_1} \sum_{s \sim S_1} \sum_{rsk \in \ci} \xi_r^*\eta_s^*\zeta_k\Phi(rsk),
  \end{equation}
  where $|\xi_r^*| \le \tau(r)^c$, $|\eta_s^*| \le \tau(s)^c$, $\max(R_1, S_1) \le P^{11/20}$,
and either $\zeta_k = 1$ for all $k$, or $|\zeta_k| \le \tau(k)^c$ and $R_1S_1 \ge P^{27/35}$.
\end{lemma}

\begin{proof}
  This can be established similarly to Lemma 5.4 in \cite{Kumchev},
which contains (essentially) the case $\Phi(m) = \chi(m)e(\beta m^2)$.
The second and third conditions in \eqref{M.13} can serve as a replacement for the hypothesis $z \le P^{23/140}$ in \cite{Kumchev}.
\end{proof}

The above result covers $\rho_2$ and $\rho_3$, while the following lemma covers additional sums that
arise in $\rho_4$.

\begin{lemma}\label{lM.1B}
Let $W \le R \le P^{1/2}$.  Then the sum
\beq\label{lb.eq}
\sum_{m \in \ci} \sum_{p \sim R}  \psi(m/p,p)
\endq
can be expressed as a linear combination of $O(\cl^c)$ sums of the form \eqref{M.14}
where the parameters satisfy the same conditions as in Lemma \ref{lM.1}.
The same conclusion is also reached for the sum
\beq\label{lc.eq}
\sum_{m \in \ci} \threelinesum{pr>W, qr>Y}{Z<r<q<Y}{p<V} \psi(m/(pqr),r).
\endq
\end{lemma}

\begin{proof} We begin with the sum \eqref{lb.eq} which clearly detects products of two primes since $p > (m/p)^{\frac12}$ here.
If $P^{9/20} \le R \le P^{1/2}$ the result is immediate with the
variable $k$ identically equal to $1$. Otherwise, let $u=P^{9/40}
R^{- 1/2}$. We apply Heath-Brown's generalized Vaughan Identity to
the variable $m/r$ as given by \cite[Lemma 2.8]{PDS} (note that
$12Nu$ there should read $N^{\frac12}u$). This gives Type II sums
with one range of size $P^{9/20}/R$ to $P^{1/3}$, and Type I sums
where the variable with an ``unknown" weight has size $\le
(P/R)^{1/2}u$. These sums are quickly shown to have the required
properties.

For the sum in \eqref{lc.eq} we immediately have a sum of the correct form when $r< P^{8/35}$ or $pqr < P^{11/20}$. For
the remainder of the sum we note that $P^{8/35} > (P/pqr)^{1/2}$ and so $\psi(m/pqr,r)$ detects primes only.
Again we can apply Heath-Brown's identity and obtain the required result.

For either of the above sums we could have used the Alternative Sieve technique we employ later, but the appeal to Heath-Brown's identity is quicker.
\end{proof}

\begin{lemma}\label{lM.2}
  Let $l \in \mathcal S_0 \cup \{ 1 \}$ and $g, n, D \in \mathbb N$.
Suppose that $\rho$ is either a convolution of the form \eqref{M.11}
 that satisfies \eqref{M.15}, or one of \eqref{lb.eq}, \eqref{lc.eq}.
Suppose also that $G, \Delta$
are reals such that $\Delta lDG^{1 + \epsilon} \le P^{-31/20}$ and that $\mathcal H(D, G)$ is
a set of characters $\chi = \xi\psi$, where $\xi$ is a character modulo $D$ and $\psi$ a primitive
character modulo $q$, with $q \le G$ and $(q, D) = 1$. Then
  \begin{equation}\label{M.16}
    \sum_{\chi \in \mathcal H(D, G)} w(q) \bigg( \int_{-\Delta}^{\Delta} |W_l(\chi, \beta)|^2 \, d\beta \bigg)^{1/2} \ll l^{-1}w(1)g^{\epsilon}\cl^c,
  \end{equation}
  where $w(q) = (n, [q, g])^{1/2}[q, g]^{-1 + \epsilon}$.
\end{lemma}

\begin{proof}
  We can use Lemma \ref{lM.0} to deduce \eqref{M.16} from the inequality
  \begin{equation}\label{M.17}
    \twolinesum{\chi \in \mathcal H(D, G)}{d \mid q} \bigg( \int_{-\Delta}^{\Delta} |W_l(\chi, \beta)|^2 \,
d\beta \bigg)^{1/2} \ll l^{-1}\big( 1 + d^{-1}lDG^2\Delta P^{31/20} \big) \cl^c.
  \end{equation}
  We first consider the case $l = 1$. If $\rho$ is of the form \eqref{M.11}, it
 follows from the hypotheses \eqref{M.15} that $R, S$ and $z(r, s)$
satisfy the hypotheses \eqref{M.13} of Lemma \ref{lM.1}. By Lemma \ref{lM.1B}
we obtain the required conclusions for the convolutions \eqref{lb.eq} and \eqref{lc.eq} too.
Thus, we may assume that $W_1(\chi, \beta)$ is
given by \eqref{M.14} with $\Phi(m) = \chi(m)e(\beta m^2)$ (with the appropriate adjustment when $\chi$ is principal).
Then the argument of Lemma 4.3 in \cite{Choi2} shows that the left side of \eqref{M.17} is bounded above by
  \begin{equation}\label{M.18}
    \Delta P\cl T^{-1} \twolinesum{\chi \in \mathcal H(D, G)}{d \mid q}
    \int_{-T}^T |F(it, \chi)| \, dt + d^{-1}G^2\Delta P^{1/2},
  \end{equation}
  where $\Delta P^2 \le T \le P^{10}$ and $F(it, \chi)$ is the Dirichlet polynomial
  \[
    F(it, \chi) = \sum_{r \sim R_1} \sum_{s \sim S_1} \sum_{k \asymp P/(R_1S_1)} \xi_r^*\eta_s^*\zeta_k \chi(rsk) (rsk)^{-it}.
  \]
  We can apply Theorem 2.1 in \cite{Choi2} to the sum in \eqref{M.18} to obtain \eqref{M.17} with $l = 1$.

  Suppose now that $l = p \in \mathcal S_0$. We note that this case cannot occur for \eqref{lb.eq}.  The following argument is for the case where
$\rho$ is of the form \eqref{M.11}, but it can easily be adapted for \eqref{lc.eq}.
We may assume that $z(r, s) \le P^{\sigma}$, since $W_p(\chi, \beta)$ is otherwise an empty sum.
The left side of \eqref{M.17} equals
  \begin{equation}\label{M.19}
    p^{-1} \twolinesum{\chi \in \mathcal H(D, G)}{d \mid q} \bigg( \int_{-\Delta_p}^{\Delta_p} |W_p(\chi, \beta p^{-2})|^2 \,
d\beta \bigg)^{1/2}, \qquad \Delta_p = p^2\Delta.
  \end{equation}
  The sum $W_p(\chi, \beta p^{-2})$ splits into three subsums: a subsum where $r = pr'$; a subsum where $p \nmid r$ and $s = ps'$;
and a subsum where $p \nmid rs$ and $k = pk'$.
Each of these three subsums can be represented in the form $W_p^*(\chi, \beta)$, where
  \beq\label{GHadd2}
    W_p^*(\chi, \beta) = \sum_{r \sim R'} \sum_{s \sim S'} \sum_{k \sim P/(prs)} \xi_r' \eta_s' \psi(rsk, z) \chi(rsk) e(\beta (rsk)^2 ),
  \endq
  with $R' \le R$, $S' \le S$, $|\xi_r'| \le \tau(r)^c$ and $|\eta_s'| \le \tau(s)^c$. We have
  \begin{align*}
    &\max(R', S') \le W \le (P/p)^{11/20} &&\text{since } 1 - 4\sigma \le \sfrac{11}{20}(1 - \sigma); \\
    &z\min(R', S') \le P^{\sigma}V \le (P/p)^{11/20} &&\text{since } 3\sigma \le \sfrac{11}{20}(1 - \sigma); \\
    &z \le P^{\sigma} \le (P/p)^{8/35} &&\text{since } \sigma \le \sfrac{8}{35}(1 - \sigma).
  \end{align*}
  We can therefore apply Lemma \ref{lM.1} to decompose $W_p^*(\chi, \beta)$ into
sums of the form \eqref{M.14} with $\Phi(m) = \chi(m)e(\beta m^2)$ and $Pp^{-1}$ in place of $P$.
By Theorem 2.1 in \cite{Choi2}, the quantity \eqref{M.19} with $W_p^*(\chi, \beta)$ in place of $W_p(\chi, \beta p^{-2})$ is bounded by
  \[
    p^{-1}\big( 1 + d^{-1}DG^2\Delta P^{31/20}p^{9/20} \big) \cl^c,
  \]
  whence \eqref{M.17} with $l = p$ follows.
\end{proof}

\subsection{Bounds for complete exponential sums}
\label{secM.2}

Given characters $\chi_1, \chi_2, \chi_3$ modulo $q$ and a vector $\mathbf b \in \mathbb Z^4$, we define
\begin{equation}\label{M.20}
  B(q, \mathbf b; \chi_1, \chi_2, \chi_3) = \frac 1{\phi(q)^3} \twolinesum{1 \le a \le q}{(a, q) = 1} S(\chi_1, ab_1) S(\chi_2, ab_2) S(\chi_3, ab_3) e_q(-ab_4).
\end{equation}
It is not difficult to express $B(q, \mathbf b; \chi_1, \chi_2, \chi_3)$ as a linear combination of Gauss sums
\[
  \tau_a(\chi) = \sum_{1 \le h \le q} \chi(h)e_q(ah).
\]
Indeed, by the orthogonality of the characters modulo $q$,
\[
  S(\chi, a) = \sum_{\xi^2 = \bar\chi} \tau_a(\xi),
\]
where the summation is over the characters $\xi$ modulo $q$ with $\xi^2 = \bar\chi$. Thus,
\begin{equation}\label{M.21}
  B(q, \mathbf b; \chi_1, \chi_2, \chi_3) = \frac 1{\phi(q)^3} \twolinesum{\xi_1, \xi_2, \xi_3}{\xi_j^2
= \bar\chi_j} \tau_{b_1}(\xi_1)\tau_{b_2}(\xi_2)\tau_{b_3}(\xi_3)\overline{\tau_{b_4}(\xi_1\xi_2\xi_3)}.
\end{equation}
We also note that $B(q)$ is multiplicative as a function of $q$ in the following sense: if $q = q_1q_2$, $(q_1, q_2)
= 1$, and $\chi_j = \chi_{j, 1} \chi_{j, 2}$ with $\chi_{j, i}$ a character modulo
$q_i$, then
\begin{equation}\label{M.22}
  B(q, \mathbf b; \chi_1, \chi_2, \chi_3) = B(q_1, \mathbf b; \chi_{1, 1}, \chi_{2, 1}, \chi_{3, 1})B(q_2, \mathbf b; \chi_{1, 2}, \chi_{2,2}, \chi_{3, 2}).
\end{equation}
The proofs of the above properties are similar to those of parts (a) and (d) of Lemma 2.5 in \cite{Choi}.
We now record upper bounds for $|B(q, \mathbf b; \chi_1, \chi_2, \chi_3)|$ for several special choices of $\mathbf b$.

\medskip

\paragraph{\em Case 1:} $\mathbf b = \mathbf b_1 = (1, 1, 1, n)$. By virtue of \eqref{M.22}, it suffices to consider the case when $q = p^e$ for some prime $p$.
We deal with the case of an odd prime $p$. The case $q = 2^e$, $e \ge 3$, can be dealt with in a similar fashion, and when $q = 2$ or $4$,
we may use the trivial bound. Using the bound (see Lemma 3.1 in \cite{LiTs05})
\[
  |\tau_a(\chi)| \le (a, q)^{1/2}q^{1/2},
\]
we deduce immediately from \eqref{M.21} that
\begin{equation}\label{M.23}
  B(p^e, \mathbf b_1; \chi_1, \chi_2, \chi_3) \le 8(n, p^e)^{1/2}p^{3 - e}\phi(p)^{-3}.
\end{equation}
In the special case when $\chi_1, \chi_2$ and $\chi_3$ are all principal, we can improve on this.
We pause at this stage to write
\[
\mathfrak s(q, n) = B(q, \mathbf b_1; \chi_0, \chi_0, \chi_0)
=  \twolinesum {1 \le a \le q}{(a, q) = 1}
  \frac {S(\chi_0, a)^3}{\phi(q)^3} e_q (-an).
\]
We also write, for future reference,
\[
\mathfrak S_3(n,Q) = \sum_{q \le Q} \mathfrak s(q,n).
\]
We then have the following result.

\begin{lemma}\label{sqn}
For all $q \ge 2$ and $n \ge 1$ we have the two estimates:
\beq\label{GHadd1}
|\mathfrak s(q,n)| \le \tau(q)^3 \frac{q^2}{\phi(q)^3},
\endq
and
\beq\label{loglog10}
|\mathfrak s(q,n)| \le \tau((q,n))^2 \frac{(\log \log q)^{10}}{q}.
\endq
Moreover, if $p^2 \mid q$ with $p > 2$, or if $16 \mid q$, we have $\mathfrak s(q,n) = 0$.
\end{lemma}

\begin{proof}
The cases $q=2^j$ can be quickly checked.  It then suffices to consider the case $q = p >2$,
since the factors $\tau_1(\xi_j)$ in \eqref{M.21} vanish when $q = p^e$, $e \ge 2$. When $q = p$ and $\chi_j$ is principal,
each $\xi_j$ is either principal or a Legendre symbol. If some $\xi_j$ is principal, we have $|\tau_1(\xi_j)| = 1$.
We also note that if $\xi$ is non-principal, then
\[
\left| \tau_n(\xi)\right| = \begin{cases}  0 \ &\text{if } p \mid n, \\  p^{1/2} \ &\text{otherwise.} \end{cases}
\]
On the other hand, if $\xi$ is principal, then
\[
\left| \tau_n(\xi)\right| = \begin{cases}  p-1 \ &\text{if } p \mid n,\\ 1  &\text{otherwise.} \end{cases}
\]
When $p \mid n$, we deduce that the modulus of the sum on the right hand side of \eqref{M.21} is
\[
  \le \begin{cases}
    p - 1     &\text{if each $\xi_j$ is principal},\\
    3p(p - 1) &\text{if exactly one $\xi_j$ is principal}, \\
    0         &\text{otherwise.} \end{cases}
\]
When $p \nmid n$, the modulus of the sum on the right hand side of \eqref{M.21} is
\[
  \le \begin{cases}
    1   &\text{if each $\xi_j$ is principal},\\
    p^2 &\text{if no $\xi_j$ is principal}, \\
    3p  &\text{otherwise.} \end{cases}
\]
We thus have
\[
  |\mathfrak s(p,n)| \le \begin{cases}
    p^2\phi(p)^{-3}(1 + 7/p) &\text{if } p \nmid n,\\
    3p^2\phi(p)^{-3}         &\text{if } p \mid n. \end{cases}
\]
The bounds \eqref{GHadd1} and \eqref{loglog10} quickly follow.
\end{proof}

We note, for future reference, that the main contribution in the case $p \nmid n$
can be explicitly calculated, namely
\beq\label{futureref}
\mathfrak s(p,n) = \left( \frac{-n}{p} \right) \frac{p^2}{(p-1)^3} + \gamma(p,n)
\endq
where $|\gamma(p,n)| \le 7p/(p-1)^3$.

Suppose now that $\chi_j$ has conductor $q_j$, and let $q_0 = [q_1, q_2, q_3]$. By \eqref{M.22}--\eqref{GHadd1},
\begin{equation}\label{M.25}
  B(q, \mathbf b_1; \chi_1, \chi_2, \chi_3) \ll (n, q_0)^{1/2}q^{-1}\tau(q)^c,
\end{equation}
whence
\begin{equation}\label{M.26}
  \twolinesum{q \le Q}{q_0 \mid q} B(q, \mathbf b_1; \chi_1\chi_0, \chi_2\chi_0, \chi_3\chi_0) \ll (n, q_0)^{1/2}q_0^{-1+\epsilon}\cl^c.
\end{equation}
Here $\chi_0$ denotes the principal character modulo $q$.

\medskip

\paragraph{\em Case 2:} $\mathbf b = \mathbf b_D = (1, p_1^2, p_2^2, n)$, where $p_1, p_2$
are distinct odd primes and $D = p_1p_2$. When $(q, D) = 1$, similarly to \eqref{M.25}, we have
\[
  B(q, \mathbf b_D; \chi_1, \chi_2, \chi_3) \ll (n, q_0)^{1/2}q^{-1}\tau(q)^c,
\]
where $q_0 = [q_1, q_2, q_3]$, $q_j$ being the conductor of $\chi_j$.
When $q = p_1$, the factor $\tau_{p_1^2}(\xi_2)$ in \eqref{M.21} vanishes
unless $\xi_2$ is principal, in which case that factor equals $\phi(p_1)$. Hence,
\begin{align*}
  |B(p_1, \mathbf b_D; \chi_1, \chi_2, \chi_3)| &\le \frac 1{\phi(p_1)^2}
  \twolinesum{\xi_1, \xi_3}{\xi_j^2 = \bar\chi_j} |\tau_{1}(\xi_1)\tau_{p_2^2}(\xi_3)\tau_{n}(\xi_1\xi_3)| \\
  &\le 4(n, p_1)^{1/2}p_1^{3/2}\phi(p_1)^{-2}.
\end{align*}
Similarly,
\[
  |B(p_2, \mathbf b_D; \chi_1, \chi_2, \chi_3)| \le 4(n, p_2)^{1/2}p_2^{3/2}\phi(p_2)^{-2}.
\]
Now, let $q = p_1p_2r$, where $(r, p_1p_2) = 1$, and suppose that $\chi_j$ has conductor $q_j$. We deduce that
\begin{equation}\label{M.27}
  B(q, \mathbf b_D; \chi_1, \chi_2, \chi_3) \ll \sqrt{D(n, D)(n, r_0)}q^{-1}\tau(q)^c,
\end{equation}
where $r_0 = ([q_1, q_2, q_3], r)$.

\medskip

\paragraph{\em Case 3:} $\mathbf b = \mathbf b_p = (1, p^2, p^2, n)$. When $q = p^e$, $e \le 2$, the factors $\tau_{p^2}(\xi_j)$, $j = 2, 3$,
in \eqref{M.21} vanish unless $\xi_j$ is principal. Hence,
\[
  |B(p^e, \mathbf b_p; \chi_1, \chi_2, \chi_3)| \le \frac 1{\phi(p^e)} \sum_{\xi_1^2 = \bar\chi_1} |\tau_{1}(\xi_1)\tau_{n}(\xi_1)| \le 2p\phi(p)^{-1},
\]
on noting that $|\tau_n(\xi_1)| \le p^{e/2}$ when $\xi_1$ is non-principal and $|\tau_1(\xi_1)| \le 1$ when $\xi_1$ is principal.
Suppose now that $q = p^er$, with $e \le 2$ and $(r, p) = 1$, and that $\chi_j$ has conductor $q_j$. Then, similarly to \eqref{M.27}, we have
\begin{equation}\label{M.28}
  B(q, \mathbf b_p; \chi_1, \chi_2, \chi_3) \ll p^e(n, r_0)^{1/2}q^{-1}\tau(q)^c,
\end{equation}
where $r_0 = ([q_1, q_2, q_3], r)$. We also remark that when $e = 2$, the left side of \eqref{M.28}
vanishes unless $p^2 \mid q_1$ and $(p, q_2q_3) = 1$, in which case $r_0 = [q_1p^{-2}, q_2, q_3]$.

\subsection{Estimation of $I_1$}
\label{secM.3}

We can rewrite $I_1$ as the multiple sum
\begin{equation}\label{M.30}
  \sum_{q \le Q} \sum_{\chi_1 \bmod q} \sum_{\chi_2 \bmod q} \sum_{\chi_3 \bmod q}
  B(q, \mathbf b_1; \chi_1, \chi_2, \chi_3) J(q, n; \chi_1, \chi_2, \chi_3),
\end{equation}
where $B(q, \mathbf b_1; \chi_1, \chi_2, \chi_3)$ is defined by \eqref{M.20} with $\mathbf b = \mathbf b_1 = (1, 1, 1, n)$ and
\[
  J(q, n; \chi_1, \chi_2, \chi_3) = \int_{-\omega/q}^{\omega/q} W_1(\chi_1, \beta)
  W_{2,1}(\chi_2, \beta) W_{4,1}(\chi_3, \beta) e(-\beta n) \, d\beta.
\]
We now pass to primitive characters in \eqref{M.30}. In general, if $\chi \bmod q$, $q \le Q$, is
induced by a primitive character $\chi^* \bmod r$, $r \mid q$, we have
\begin{equation}\label{M.31}
  W_1(\chi, \beta) = W_1(\chi^*, \beta)
\end{equation}
and
\begin{equation}\label{M.32}
  W_{j,1}(\chi, \beta) \ll |W_{j,1}(\chi^*, \beta)| + \twolinesum{p \in \mathcal S_q}{p \nmid r}|W_{j,p}(\chi^*, \beta)| + E(q, r),
\end{equation}
where $E(q, r)$ denotes the number of integers $m \in \ci$ with $(m, r) = 1$, $\psi(m, Z) = 1$
and $(m, q) \ge P^{\sigma}$. Since $Z^3 > Q$, if an integer $m$ is counted in $E(q, r)$, then $(m, q)$ is either a prime $p \ge P^{\sigma}$ or the product $p_1p_2$ of two distinct primes $p_1, p_2 \ge Z$. Now, given a character $\chi$ modulo $r$, we define
\begin{gather*}
  W_0(\chi) = \max_{|\beta|\le \omega/r} |W_1(\chi, \beta)|, \quad W_{j, l}(\chi) = \bigg( \int_{-\omega/lr}^{\omega/lr} |W_{j,l}(\chi, \beta)|^2 d\beta \bigg)^{1/2}, \\
  W_j(\chi) = W_{j, 1}(\chi), \quad W_j^{\sharp}(\chi) = \sum_{p \in \mathcal S_0} W_{j, p}(\chi).
\end{gather*}
Let $\chi_j^*$ denote the primitive character modulo $q_j$, $q_j \mid q$, inducing $\chi_j$. By \eqref{M.31} and \eqref{M.32},
\[
  J(q, n; \chi_1, \chi_2, \chi_3) \ll \sum_{1 \le i \le 9} J_i(q; \chi_1^*, \chi_2^*, \chi_3^*),
\]
where each $J_i(q; \chi_1^*, \chi_2^*, \chi_3^*)$ is a product of the form $W_0(\chi_1^*)W_2^{\flat}(\chi_2^*)W_4^{\flat}(\chi_3^*)$,
with $W_j^{\flat}(\chi)$ one of the following:
\[
  W_j(\chi), \quad W_j^{\sharp}(\chi), \quad (\omega/q)^{1/2}E(q, q_j).
\]

Suppose first that $J_i(q; \chi_1^*, \chi_2^*, \chi_3^*)$ is one of
the four products involving only $W_j(\chi)$ and
$W_j^{\sharp}(\chi)$. We note that in this case $J_i(q; \chi_1^*,
\chi_2^*, \chi_3^*)$ depends only on the characters and not on $q$.
Thus, its contribution to the final bound for \eqref{M.30} is
bounded above by
\begin{gather}
  \sideset{}{^*}\sum_{q_1, \chi_1} \sideset{}{^*}\sum_{q_2, \chi_2} \sideset{}{^*}\sum_{q_3, \chi_3}
  J_i(\chi_1, \chi_2, \chi_3) B_1(\chi_1, \chi_2, \chi_3), \label{M.33}
\end{gather}
where $\sum_{q_j, \chi_j}^*$ denotes a summation over the primitive characters of moduli $q_j \le Q$,
and $B_1(\chi_1, \chi_2, \chi_3)$ is the sum in \eqref{M.26}. Hence, by \eqref{M.26}, the sum \eqref{M.33} is bounded by
\begin{equation}\label{M.34}
  \cl^c \sideset{}{^*}\sum_{q_1, \chi_1} \sideset{}{^*}\sum_{q_2, \chi_2} \sideset{}{^*}\sum_{q_3, \chi_3}
  (n, q_0)^{1/2}q_0^{-1 + \epsilon} J_i(\chi_1, \chi_2, \chi_3),
\end{equation}
where $q_0 = [q_1, q_2, q_3]$. The four such sums can be estimated in a similar fashion, so we present only the details of the estimation of
\[
  \sideset{}{^*}\sum_{q_1, \chi_1} \sideset{}{^*}\sum_{q_2, \chi_2} \sideset{}{^*}\sum_{q_3, \chi_3}
  (n, q_0)^{1/2}q_0^{-1 + \epsilon} W_0(\chi_1) W_2^{\sharp}(\chi_2) W_4(\chi_3).
\]
Since $\rho_4$ satisfies hypothesis (iv), Lemma \ref{lM.2} with $l = 1$ gives
\[
  \sideset{}{^*}\sum_{q_3, \chi_3} (n, q_0)^{1/2}q_0^{-1 + \epsilon} W_4(\chi_3) \ll (n, \tilde q_0)^{1/2}\tilde q_0^{-1 + 2\epsilon}\cl^c,
\]
where $\tilde q_0 = [q_1, q_2]$. Furthermore, since $\rho_2$ satisfies hypothesis (iv), Lemma \ref{lM.2} with $l = p$, $p \in \mathcal S_0$, gives
\begin{equation}\label{M.35}
  \sideset{}{^*}\sum_{q_2, \chi_2} (n, \tilde q_0)^{1/2}\tilde q_0^{-1 + 2\epsilon} W_2^{\sharp}(\chi_2)
  \ll (n, q_1)^{1/2}q_1^{-1 + 3\epsilon}\cl^c.
\end{equation}
Finally, by Lemma 2.3 in \cite{Zhan2},
\begin{equation}\label{M.36}
  \sideset{}{^*}\sum_{q_1, \chi_1} (q_1, n)^{1/2}q_1^{-1+3\epsilon} W_0(\chi) \ll P\cl^{-A}
\end{equation}
for any fixed $A > 0$.

Next, we estimate the contribution to \eqref{M.30} from a product $J_i(q; \chi_1^*, \chi_2^*, \chi_3^*)$ where at
 least one of the factors $W_j^{\flat}(\chi)$ is of the form $(\omega/q)^{1/2}E(q, q_j)$. Let us consider, for example, the contribution from the product
\[
  W_0(\chi_1^*)W_2(\chi_2^*)(\omega/q)^{1/2}E(q, q_3).
\]
By \eqref{M.25}, this contribution does not exceed
\begin{equation}\label{M.37}
  \omega^{1/2} \sideset{}{^*}\sum_{q_1, \chi_1} \sideset{}{^*}\sum_{q_2, \chi_2} \sideset{}{^*}\sum_{q_3, \chi_3}
  W_0(\chi_1)W_2(\chi_2)(n, q_0)^{1/2}\twolinesum{q \le Q}{q_0 \mid q} E(q, q_3)q^{-3/2 + \epsilon},
\end{equation}
where $q_0 = [q_1, q_2, q_3]$. Let $\mathcal D$ denote the set of integers $d \le Q$ that
are either a prime $p \ge P^{\sigma}$ or a product $p_1p_2$ of two primes $p_1, p_2 \ge Z$. The innermost sum in \eqref{M.37} is bounded by
\begin{align*}
  \sum_{d \in \mathcal D} \frac Pd \twolinesum{q \le Q}{[q_0, d] \mid q} q^{-3/2 + \epsilon}
  \ll P \twolinesum{d \in \mathcal D}{dq_3 \le Q} d^{-1} [q_0, d]^{-3/2 + \epsilon} = \Sigma(q_3), \quad \text{say}.
\end{align*}
Put $\tilde q_0 = [q_1, q_2]$. Summing this bound over $q_3$, we find that
\begin{align}
  \sideset{}{^*}\sum_{q_3, \chi_3} (n, q_0)^{1/2} \Sigma(q_3) &\ll
  P \sum_{d \in \mathcal D} \frac 1d \sum_{q_3 \le Q/d} \frac {(n, q_0)^{1/2}q_3}{[q_0, d]^{3/2 - \epsilon}} \notag \\
  &\ll Q^{1/2}P \sum_{d \in \mathcal D} \frac 1{d^2} \sum_{q_3 \le Q/d} \frac {(n, q_0)^{1/2}}{[q_0, d]^{1 - \epsilon}}, \notag \\
  &\ll Q^{1/2}P \sum_{d \in \mathcal D}  \frac {(\tilde q_0, d)^{1 - \epsilon}}{d^{3 - \epsilon}} \sum_{q_3 \le Q/d} (n, q_0)^{1/2}q_0^{-1 + \epsilon} \notag \\
  &\ll Q^{1/2 + \epsilon}P (n, \tilde q_0)^{1/2}\tilde q_0^{-1 + 2\epsilon} \sum_{d \in \mathcal D} \frac{(\tilde q_0, d)^{1 - \epsilon}}{d^3} \notag \\
  &\ll Q^{1/2 + \epsilon}P^{1 - 2\sigma} (n, \tilde q_0)^{1/2}\tilde q_0^{-1 + 2\epsilon}, \label{M.38}
\end{align}
where we have used that $q_3[q_0, d]^{-1/2} \le Q^{1/2}d^{-1}$ and $(q_0, d) = (\tilde q_0, d)$.
Combining \eqref{M.38}, the variant of \eqref{M.35} for $W_2(\chi_2)$, and \eqref{M.36}, we conclude that the quantity \eqref{M.37} does not exceed
\[
  \omega^{1/2}P^{2 -2\sigma + \epsilon}Q^{1/2} \ll QP^{1 - 2\sigma + 2\epsilon} \ll P^{1 - \epsilon}.
\]
This completes the estimation of $I_1$.

\subsection{Estimation of $I_2$.}
\label{secM.4}

We first consider the part of $I_2$ where $p_1 \ne p_2$. For every such pair of primes $\mathbf p = (p_1, p_2)$,
the integral over $\mathfrak M_{\mathbf p}$ equals
\begin{equation} \label{M.40}
  \twolinesum{q \le Q}{D \mid q}
  \sum_{\chi_1 \bmod q} \sum_{\chi_2 \bmod q/{p_1}} \sum_{\chi_3 \bmod q/{p_2}} B(q, \mathbf b_D; \chi_1, \chi_2, \chi_3)J(q, n, D; \chi_1, \chi_2, \chi_3),
\end{equation}
where $D = p_1p_2$,
\begin{gather*}
  B(q, \mathbf b_D; \chi_1, \chi_2, \chi_3) = \frac {\phi(D)}{\phi(q)^3} \twolinesum{1 \le a \le q}{(a, q) = 1}  S(\chi_1, a) S(\chi_2, ap_1) S(\chi_3, ap_2) e_q(-an), \\
  J(q, n, D; \chi_1, \chi_2, \chi_3) =  \int_{-\omega/q}^{\omega/q} W_1(\chi_1, \beta)
  W_{2,p_1}(\chi_2, \beta) W_{4,p_2}(\chi_3, \beta) e(-\beta n) \, d\beta.
\end{gather*}
When $p \mid q$ and $\chi$ is a character modulo $q_p$, we have
\begin{equation}\label{M.41}
  S(\chi, ap) = \begin{cases}
    \phi(p)^{-1} S(\chi\chi_0, ap^2) & \text{if } p^2 \nmid q, \\
    p^{-1}S(\chi\chi_0, ap^2) & \text{if } p^2 \mid q,
  \end{cases}
\end{equation}
where $\chi\chi_0$ is the character modulo $q$ induced by $\chi$.
Therefore, $B(q, \mathbf b_D; \chi_1, \chi_2, \chi_3)$ is, in fact, the sum \eqref{M.20} with $\mathbf b = \mathbf b_D = (1, p_1^2, p_2^2, n)$.

Similarly to \eqref{M.32}, we have
\begin{equation}\label{M.43}
  W_{j,p}(\chi, \beta) \ll |W_{j,p}(\chi^*, \beta)| + E(q, r)
\end{equation}
whenever $\chi$ is character modulo $q_p$ induced by a character $\chi^*$ modulo $r$.
Using \eqref{M.31}, \eqref{M.43} and \eqref{M.27}, we reduce the estimation of \eqref{M.40} to the estimation of four sums of the form
\begin{equation}\label{M.42}
  \sqrt{D(n,D)} \sideset{}{^*}\sum_{q_1, \chi_1} \sideset{}{^*}\sum_{q_2, \chi_2} \sideset{}{^*}\sum_{q_3, \chi_3}
  (n, r_0)^{1/2}W_0(\chi_1) \twolinesum{q \le Q}{q_0 \mid q} W_{2,p_1}^{\flat}(\chi_2)W_{4,p_2}^{\flat}(\chi_3)q^{-1 + \epsilon},
\end{equation}
where $q_0 = [q_1, p_1q_2, p_2q_3]$, $r_0 = q_0/D$, $W_{2,p_1}^{\flat}(\chi_2)$ represents
either $W_{2,p_1}(\chi_2)$ or $(\omega/q)^{1/2}E(q, q_2)$, and $W_{4,p_2}^{\flat}(\chi_3)$ is defined similarly to $W_{2,p_1}^{\flat}(\chi_2)$.
The contribution from the sum involving the factor $W_{2,p_1}(\chi_2)W_{4,p_2}(\chi_3)$ is bounded by
\begin{equation}\label{M.44}
  Q^{\epsilon} \sqrt{\frac {(n, D)}D} \sideset{}{^*}\sum_{q_1, \chi_1} \sideset{}{^*}\sum_{q_2, \chi_2} \sideset{}{^*}\sum_{\substack{q_3, \chi_3\\ q_0 \le Q}}
  (n, r_0)^{1/2}r_0^{-1} W_0(\chi_1) W_{2,p_1}(\chi_2)W_{4,p_2}(\chi_3).
\end{equation}
The condition $[q_1, p_1q_2, p_2q_1] \le Q$ implies that each character $\chi_3$ is either primitive with a modulus $r_3 \le QD^{-1}$,
or the product of such a character and a (primitive) character modulo $p_1$. Thus, the sum over $q_3$ and $\chi_3$ splits into
two sums of the form appearing in Lemma \ref{lM.2}: one over $\chi_3 \in \mathcal H(1, QD^{-1})$ and one over $\chi_3 \in \mathcal H(p_1, QD^{-1})$.
Similarly, the sum over $q_2$ and $\chi_2$ splits into sums over $\chi_2 \in \mathcal H(1, QD^{-1})$ and $\chi_2 \in \mathcal H(p_2, QD^{-1})$,
and the sum over $q_1$ and $\chi_1$ splits into four sums over the sets $\mathcal H(d, QD^{-1})$, $d \mid D$. Observe that $r_0 = [r_1, r_2, r_3]$,
where the $r_j$'s are the moduli of the primitive characters with moduli $\le QD^{-1}$.
Hence, Lemma \ref{lM.2} with $l = p_1$, $D = l \in \{1, p_1\}$ and $\Delta = \omega/(p_2r_3l)$ gives
\[
  \sum_{\chi_3 \in \mathcal H(l, QD^{-1})} (n, r_0)^{1/2}r_0^{-1} W_{4,p_2}(\chi_3) \ll p_2^{-1}(n, [r_1, r_2])^{1/2}[r_1, r_2]^{-1 + \epsilon} \cl^c.
\]
Another application of Lemma \ref{lM.2} to the sum over $\chi_2$ and an application of a variant of Lemma 4.1 in \cite{Choi2} (a combination of that lemma with Lemma \ref{lM.0} above) to the sum over $\chi_1$ show that that the sum \eqref{M.44} is $\ll (n, D)^{1/2}D^{-3/2}QP^{1 + \epsilon}$.

As to the estimation of the remaining three sums of the form \eqref{M.42}, we note that the condition $q_0 \mid q$ implies that
$E(q, q_j) < PD^{-1}$  and $q_j \le QD^{-1}$.
Hence, the contribution from the product involving $E(q, q_2)E(q, q_3)$, for example, is bounded above by
\begin{align*}
   \omega P^2\sqrt{\frac {(n, D)}{D^3}} \sideset{}{^*}\sum_{q_1, \chi_1} W_0(\chi_1) \sum_{q_2, q_3 \le QD^{-1}} q_2q_3 (n, r_0)^{1/2} q_0^{-2 + \epsilon} \\
  \ll  \omega QP^{2 + \epsilon}\sqrt{\frac {(n, D)}{D^9}}\sum_{l \mid D} \sum_{\chi \in \mathcal H(l, QD^{-1})} (n, r)^{1/2}r^{-1 + 3\epsilon} W_0(\chi).
\end{align*}
Another appeal to the variant of Lemma 4.1 in \cite{Choi2} used above shows that the last expression is
\[
  \ll \omega QP^{3 + 2\epsilon}\sqrt{\frac {(n, D)}{D^9}} \ll (n, D)^{1/2}D^{-3/2}Q^2P^{1 + 3\epsilon}Z^{-6}.
\]
Therefore, the total contribution to $I_2$ from pairs $(p_1, p_2)$ of distinct primes is
\[
  \ll P^{1 + 3\epsilon} \sum_{p_1, p_2 \in \mathcal S_0} (n, p_1p_2)^{1/2}(p_1p_2)^{-3/2} \ll P^{1 + 3\epsilon}Z^{-1}.
\]

Finally, let $p \in \mathcal S_0$ and $\mathbf p = (p, p)$. Then the integral over $\mathfrak M_{\mathbf p}$ appearing in $I_2$ can be expressed as
\begin{equation}\label{M.45}
  \twolinesum{q \le Q}{p \mid q}
  \sum_{\chi_1 \bmod q} \sum_{\chi_2 \bmod q_p} \sum_{\chi_3 \bmod q_p} B(q, \mathbf b_p; \chi_1, \chi_2, \chi_3)J(q, n, p; \chi_1, \chi_2, \chi_3),
\end{equation}
where $q_p = qp^{-1}$,
\begin{gather*}
  B(q, \mathbf b_p; \chi_1, \chi_2, \chi_3) = \frac 1{\phi(q)\phi(q_p)^2} \twolinesum{1 \le a \le q}{(a, q) = 1}  S(\chi_1, a) S(\chi_2, ap) S(\chi_3, ap) e_q(-an), \\
  J(q, n, p; \chi_1, \chi_2, \chi_3) =  \int_{-\omega/q}^{\omega/q} W_1(\chi_1, \beta)
  W_{2,p}(\chi_2, \beta) W_{4,p}(\chi_3, \beta) e(-\beta n) \, d\beta.
\end{gather*}
By \eqref{M.41}, $B(q, \mathbf b_p; \chi_1, \chi_2, \chi_3)$ is the sum \eqref{M.20} with $\mathbf b = \mathbf b_p = (1, p^2, p^2, n)$.
Hence, by \eqref{M.31}, \eqref{M.43} and \eqref{M.28}, the contribution to \eqref{M.45} from moduli $q$ divisible by $p$ but not
by $p^2$ does not exceed the linear combination of four sums of the form
\begin{equation}\label{M.46}
  p \sideset{}{^*}\sum_{q_1, \chi_1} \sideset{}{^*}\sum_{q_2, \chi_2} \sideset{}{^*}\sum_{q_3, \chi_3}
  (n, r_0)^{1/2}W_0(\chi_1) \twolinesum{q \le Q}{q_0 \mid q, p^2 \nmid q} W_{2,p}^{\flat}(\chi_2)W_{4,p}^{\flat}(\chi_3)q^{-1 + \epsilon},
\end{equation}
where $q_0 = [q_1, pq_2, pq_3]$, $r_0 = q_0p^{-1}$, and $W_{j,p}^{\flat}(\chi)$ has the same meaning as in \eqref{M.42}. The sum \eqref{M.46}
involving the product $W_{2,p}(\chi_2)W_{4,p}(\chi_3)$ is bounded by
\[
  Q^{\epsilon} \sideset{}{^*}\sum_{q_1, \chi_1} \sideset{}{^*}\sum_{q_2, \chi_2} \sideset{}{^*}\sum_{\substack{q_3, \chi_3\\ q_0 \le Q, p^2 \nmid q_0}}
  (n, r_0)^{1/2}r_0^{-1} W_0(\chi_1)W_{2,p}(\chi_2)W_{4,p}(\chi_3).
\]
The conditions $q_0 \le Q$ and $p^2 \nmid q$ imply that $q_2, q_3 \le Qp^{-1}$ and that $q_1 = r_1$ or $q_1 = pr_1$, where $(p, r_1) = 1$.
Furthermore, $r_0 = [r_1, q_2, q_3]$. Thus, we can again use Lemma \ref{lM.2}
and a variant of Lemma 4.1 in \cite{Choi2} to show that the last sum is $\ll p^{-2}P^{1 + \epsilon}$. The sums \eqref{M.46}
involving factors $(\omega/q)E(q, q_j)$ satisfy the same bound. For example, one of those does not exceed
\begin{align*}
  \omega^{1/2}P \sideset{}{^*}\sum_{q_1, \chi_1} \sideset{}{^*}\sum_{q_2, \chi_2} W_0(\chi_1)W_{2,p}(\chi_2) \sum_{p_1 \ge Z}
\frac 1{p_1} \sum_{q_3 \le Q/(pp_1)} (n, r_0)^{1/2}q_3 \twolinesum{q \le Q}{s \mid q} q^{-3/2 + \epsilon} \\
  \ll \omega^{1/2}P \sideset{}{^*}\sum_{q_1, \chi_1} \sideset{}{^*}\sum_{\substack{q_2, \chi_2\\ [q_1, pq_2] \le Q}}
W_0(\chi_1)W_{2,p}(\chi_2) \sum_{p_1 \ge Z} \frac {Q^{1/2}}{pp_1^2} \sum_{q_3 \le Q/(pp_1)} (n, r_0)^{1/2}s^{-1 + \epsilon},
\end{align*}
where $s = [q_1, pq_2, pp_1q_3]$. Another application of Lemma \ref{lM.2} and of the same variant of Lemma 4.1 in \cite{Choi2} as before show
that the last sum is $\ll p^{-2}P^{1 + \epsilon}QZ^{-3}$. Therefore, the total contribution to \eqref{M.45} from moduli not
divisible by $p^2$ is $\ll p^{-2}P^{1 + \epsilon}$.

Finally, we consider the contribution to \eqref{M.45} from moduli $q$ divisible by $p^2$. For such moduli,
the term $E(q, r)$ in \eqref{M.43} is superfluous. Thus, by \eqref{M.28} and the remark following it, this contribution is bounded by
\[
  p^2Q^{\epsilon} \sideset{}{^*}\sum_{q_1, \chi_1} \sideset{}{^*} \sum_{q_2, \chi_2} \sideset{}{^*}
\sum_{ \substack{ q_3, \chi_3\\ q_0 \le Q}} (n, r_0)^{1/2}q_0^{-1 + \epsilon} W_0(\chi_1)W_{2,p}(\chi_2)W_{4,p}(\chi_3)
\]
where $q_0 = [q_1, pq_2, pq_3]$, $r_0 = q_0p^{-2}$ and the moduli $q_1, q_2, q_3$ satisfy the conditions $p^2 \mid q_1$, $(p, q_2q_3) = 1$.
We note that these conditions imply that $q_2, q_3 \le Qp^{-2}$ and that $\chi_1 \in \mathcal H(p^2, Qp^{-2})$.
Thus, once again, we can use Lemma \ref{lM.2} and a variant of Lemma 4.1 in \cite{Choi2} to show
that the last sum does not exceed $\ll p^{-2}P^{1 + \epsilon}$. Therefore,
the total contribution to $I_2$ from pairs $(p, p)$, with $p \in \mathcal S_0$ is
\[
  \ll P^{1 + \epsilon} \sum_{p \in \mathcal S_0} p^{-2} \ll P^{1 + \epsilon}Z^{-1}.
\]

We have thus shown that the integral in \eqref{M.3} is $O(PL^{-A})$.  We then have
\[
\int_{\fM} f_1^*(\alpha)f_2^*(\alpha)f_4^*(\alpha) e(-\alpha n) \, d\alpha
=
\sum_{q \le Q} \mathfrak s(q,n) I(q,n),
\]
where
\[
I(q,n) =
\int_{-\omega/q}^{\omega/q} f_1^*(\beta)f_2^*(\beta)f_4^*(\beta) e(-\beta n) \, d\beta.
\]
We can then use standard major arc techniques to show that, for $|\beta| < \omega/q$, we have
\[\begin{split}
f_j^*(\beta) &= \sum_{\sqrt u \in \ci} \frac{\delta_j(\sqrt{u})}{2 \sqrt{u}} e(\beta u) + O(P \omega/q)\\
&\ll
\frac{P}{P^2|\beta| + 1} + O(P \omega/q).
\end{split}
\]
The error arising from any terms involving $O(P \omega/q)$ will be smaller than the other errors which arise.
We then complete the integral over $[-\omega/q, \omega/q]$ to an integral over $[-1/2,1/2]$
incurring an error bounded by a constant times
\[
\sum_{q \le Q} |\mathfrak s(q,n)| \frac{1}{P^3} \frac{q^2}{\omega^2}
\ll \sum_{q \le Q} \tau(q)^3 \frac{q^2}{\phi(q)^3} q^2 Q^{-2} P^{1-2\epsilon}
\ll P^{1 -  \epsilon}
\]
using \eqref{GHadd1}.  This shows the main term to be
\[
\sum_{q \le Q} \mathfrak s(q,n) \twolinesum{\sqrt{m_j} \in \ci}{m_1+m_2+m_3 = n}
\frac {\delta_1(\sqrt{m_1}) \delta_2(\sqrt{m_2}) \delta_4(\sqrt{m_3})}{8 \sqrt{m_1 m_2 m_3}}
\]
Clearly we can write the main term as
\[
\mathfrak S_3(n, Q) C_2 C_4 P L^{-3} K_n (1 + o(1)).
\]
where $1 \ll K_n \ll 1$ with absolute constants.
As indicated earlier, similar but simpler working leads to an analogous result for
\[
  \int_{\fM} f_1(\alpha)^2g_3(\alpha)e(-\alpha n) \, d\alpha,
\]
with a main term $\mathfrak S_3(n, Q) C_3 P L^{-3} K_n (1 + o(1))$.
Thus we obtain
\begin{multline*}
  \int_{\mathfrak M} \left(f_1^2(\alpha) g_3(\alpha) -  f_1 (\alpha) g_2(\alpha) g_4(\alpha)\right) \,d\alpha \\
  = \mathfrak S_3(n, Q)(C_3 - C_2C_4)K_n PL^{-3}(1 + o(1)) + O(PL^{-A}).
\end{multline*}

\subsection{The singular series}
\label{secM.5}
Our goal in this section is to prove the following result.

\begin{lemma}\label{MAlem}
Write
\[
G(\alpha) = f_1(\alpha)\left(f_1(\alpha) g_3(\alpha) - g_2(\alpha)g_4(\alpha)\right).
\]
Then, for all but $O(N^{1-\sigma + \epsilon})$ integers $n \in \mathcal B$, we have
\beq\label{MA}
\int_{\fM} G(\alpha) e(-\alpha n)  \, d\alpha \gg (C_3 - C_2 C_4)P L^{-6}.
\endq
\end{lemma}

Write
\[
\Pi(n,Q) = \begin{cases} 8 \prod\limits_{3 \le p \le Q} \left(1 + \mathfrak s(p,n)  \right) \ &\text{if $n \in \ca_3$,}\\
0 &\text{otherwise.}
\end{cases}
\]
Lemma \ref{MAlem} will thus follow from our previous work once we demonstrate the following.

\begin{lemma}
For all but $O(N^{1 + \epsilon} Q^{-1})$ integers in $\cb$ we have
\beq\label{SNP}
\mathfrak S(n,Q) = \Pi(n,Q) + O\left(\exp\left(- (\log L)^{1 + \epsilon}  \right)  \right).
\endq
\end{lemma}

\begin{remark}  The reader will note in the proof that the value $Q$ can be taken as large as
$N^{1/5}$ in this part of argument.
\end{remark}

\begin{proof}
In the following we can assume that whenever the variable $q$ appears it has no square odd factor exceeding $1$ and
is not divisible by $16$.  We write
\[
\Psi(r,z) = \begin{cases}
1 &\text{if $p|r \Rightarrow p \le z$,}\\ 0 &\text{otherwise.}
\end{cases}
\]
Let $R$ be a parameter exceeding $Q$ to be determined later.  We begin by writing
\[
\begin{split}
\Pi(n,Q) - \mathfrak S(n,Q)
&=
\sum_{q > Q} \mathfrak s(q,n) \Psi(q,Q)\\
&=
\Sigma_1(n) + \Sigma_2(n)
\end{split}
\]
where $Q < q \le R$ in $\Sigma_1(n)$ and $q > R$ in $\Sigma_2(n)$.  We now use \eqref{loglog10}
to obtain
\[
\Sigma_2(n) \ll \sum_{d|n} \mu^2(d) \tau^2(d) \twolinesum{q>R}{d|q}
\Psi(q,Q) \frac{(\log \log q)^{10}}{q}.
\]
From our restriction on $q$ we note that $\Psi(q,Q)$ vanishes when $q \ge \exp(2Q)$.  Hence
\[
\Sigma_2(n) \ll L^{10} \sum_{d|n} \frac{\tau^2(d)}{d} \sum_{q > R/d} \frac{\Psi(q,Q)}{q}.
\]
We now choose $R$ to satisfy
\[
\log R = L (\log L)^{1 + 2 \epsilon}.
\]
Then, using standard bounds on the number of integers up to $2^j$ having all their
prime factors $\le Q$ from \cite{HT}, we have
\[
\begin{split}
\Sigma_2(n) &\ll L^{14} \sum_{q > R/n} \frac{\Psi(q,Q)}{q} \ll
L^{14} \sum_{2^j > R/n} 2^{-j} \sum_{q \le 2^j} \Psi(q,Q)\\
&\ll
L^{14} \sum_{2^j > R/n} \exp\left(- \frac{\log 2^j}{\log Q}  \right) \ll
\exp\left(-(\log L)^{1 + \epsilon} \right).
\end{split}
\]
We thus have
\[
\Pi(n,Q) - \mathfrak S(n,Q) = \Sigma_1(n) + O\left(\exp\left(-(\log L)^{1 + \epsilon} \right)\right).
\]
Recalling \eqref{futureref} we define $\gamma(p,n)$ for $p > 2$ by
\[
  \gamma(p,n) = \mathfrak s(p,n) - \left(\frac{-n}{p}\right) \frac{p^2}{(p-1)^3},
\]
and extend this definition to obtain a multiplicative function $\gamma(q,n)$ defined on
odd square-free $q$.

We have
\[
\sum_{Q < q \le R} \Psi(q,Q) \mathfrak s(q, n) \ll \left|\sum_{Q' < q \le R'} \theta_q \mathfrak s(q,n)  \right|,
\]
where $\theta_q \in \{0,1\}$ and $\theta_q = 0$ unless $q$ is odd and square-free.  Also
$2^j (Q',R') = (Q,R)$ for some $j \in \{0,1,2,3\}$. For the values of $q$ of relevance we have
\[
\mathfrak s(q,n) = \sum_{d|q} \frac{d^2}{\phi(d)^3} \left(\frac{-n}{d}\right) \gamma(q/d,n).
\]
It then suffices to estimate
\[
\Sigma_3(n) = \sum_{Q < dq \le R} \theta_q \theta_d \frac{q^2}{\phi(q)^3} \left(\frac{-n}{q}\right)\gamma(d,n).\]
We write $\Sigma_4(n)$ to be the subsum of this expression with $d > N^{\epsilon}$.  Then
\[
\Sigma_4(n) \ll N^{-\epsilon} \sum_{dq} d \theta_q \theta_d \frac{q^2}{\phi(q)^3} \left(\frac{-n}{d}\right) |\gamma(d,n)|.
\]
It then follows that
\[
\Sigma_4(n) \ll N^{- \epsilon} \prod_{3 \le p \le Q} \left(1 + \frac{p^2}{(p-1)^3}  \right)
\prod _{3 \le p \le Q} \left(1 + p |\gamma(p,n)|  \right).
\]
From our earlier work we know that
\[
  |\gamma(p,n)| \le \begin{cases}  3p^2(p-1)^{-3} &\text{if $p \mid n$,}\\
  7p^{-2} & \text{if $p \nmid n$.}
\end{cases}
\]
Hence,
\[
\begin{split}
\Sigma_4(n) &\ll N^{-\epsilon} \tau(n)^2 \prod_{p \le Q} \left(1+\frac{1}{p}\right) \prod_{p \le Q}\left(1 + \frac{7}{p}  \right)\\
&\ll N^{-\epsilon} \tau(n)^2 L^8 \ll N^{-\epsilon/2}.
\end{split}
\]
It therefore remains to bound
\[
\Sigma_5(n) = \sum_{d \le N^{\epsilon}} |B(d,n)| \left|\sum_{Q < qd \le R} \theta_q \frac{q^2}{\phi(q)^3} \left(\frac{-n}{d}\right)   \right|.
\]
Let $\Sigma_5(n,G) $ denote the part of the inner sum with $G < q \le 2G$.  So
\[
\Sigma_5(n) \ll \cl \sum_{d \le N^{\epsilon}} \frac{\tau((d,n))^2}{d} \sum_G \left|\Sigma_5(G,n)\right|.
\]
Write
\[
\ce(G) = \left\{n \in \ca_3: n \le N, |\Sigma_5(G,n)| \ge N^{-\epsilon/3} \right\}.
\]
We will obtain two different bounds for the cardinality of this set: one to cover the range $G \le N/Q$, the other
for the remaining values of $G$.

First consider the sum
\[
\begin{split}
\sum_{n \le N} \left|\Sigma_5(G,n)  \right|^2
&=
\sum_{n \le N}\left|\sum_{q \sim G} \theta_q \frac{q^2}{\phi(q)^3} \left(\frac{-n}{q}\right)  \right|^2\\
&=
\sum_{q_j \sim G} \theta_{q_1} \theta_{q_2} \frac{(q_1 q_2)^2}{\phi(q_1)^3 \phi(q_2)^3}
\sum_{n \le N}\left(\frac{-n}{q_1 q_2}\right).
\end{split}
\]
The terms with $q_1 = q_2$ can only be estimated trivially.  They give a contribution
\[
\le \sum_{q \sim G} \frac{q^4}{\phi(q)^6} N \ll \frac{N}{G}L.
\]
On the other hand, if $q_1 \ne q_2$ then $\big(\frac{-n}{q_1 q_2}\big)$ is a non-principal character
$\modulo {q_1q_2}$.  Hence, these terms contribute
\[
\ll \sum_{q_j \sim G} (q_1 q_2)^{\epsilon - 1/2} \ll G^{1 + \epsilon}.
\]
We thus have
\[
\sum_{n \le N} \left|\Sigma_5(G,n)  \right|^2 \ll (N/Q + G)N^{\epsilon}.
\]
It follows that
\[
|\ce(G)| \ll N^{1 + \epsilon}/Q
\]
if $Q \le G \le N/Q$.

Now we use the method of \cite{Leung} to estimate
\[
\sum_{n \le N} \left|\Sigma_5(G,n)  \right|.
\]
In the following the parameter $m$ will satisfy $2 \le m \ll (\log L)^c$.
From \cite[Lemma 6.5]{Leung} we have
\[
\begin{split}
&\sum_{n \le N} \left|\Sigma_5(G,n)\right|\\
&\ll
\left(N + G^{\frac{1}{m}}N^{\frac12}\right)(m \log N + 1)^{\frac{m^2-1}{2m}}
\left(\sum_{q \sim G} \left(\frac{q^2}{\phi(q)^3}\right)^{\frac{2m}{2m-1}}\right)^{\frac{2m-1}{2m}}\\
&\ll
\left(N + G^{\frac{1}{m}}N^{\frac12}\right) L^m (\log \log G)^3
\left(\sum_{q \sim G} q^{-\frac{2m}{2m-1}}\right)^{\frac{2m-1}{2m}}\\
&\ll
\left(N G^{-\frac{1}{2m}} + G^{\frac{1}{2m}}N^{\frac12}\right) L^{m+1} \ll
\left(N G^{-\frac{1}{2m}} + G^{\frac{1}{2m}}N^{\frac12}\right)N^{\epsilon}.
\end{split}
\]
Hence, we obtain
\[
|\ce(G)| \ll N^{1 + \delta} G^{-\frac{1}{2m}} + G^{\frac{1}{2m}}N^{\frac12 + \epsilon} \ll N^{1 + \epsilon}/Q,
\]
provided that $Q^{2m} \le G \le N^mQ^{-2m}$. Thus, so long as $Q \le N^{1/5}$, this covers the whole range from
$Q^2$ to $R$.

We can now combine our two bounds to obtain
\[
|\ce(G)| \ll N^{1 + \epsilon}Q^{-1}
\]
for each of the $O(\cl^2)$ choices for $G$ in the range $Q \le G \le R$.  This gives
\[
\sum_G \left| \Sigma_5(G,n)  \right| \ll \cl^2 N^{- \epsilon/3}
\]
for all but $O(N^{1 + \delta}Q^{-1})$ integers $n \in \cb$.
It follows that
\[
\Sigma_1(n) \ll N^{-\epsilon/4}
\]
for all but $O(N^{1 + \epsilon}Q^{-1})$ integers in $n \in \cb$,
which completes the proof.
\end{proof}

\section{Minor arc estimates for $f_j(\alpha)$}

One limit to the exponent saving that we can make in our theorems comes from
the best bounds we can obtain for one of the $f_j(\alpha)$ on the minor arcs.
The best result at the moment for $f_1(\alpha)$ would lead to only
$\sigma = 1/8$.  The combination of a sieve method and bilinear exponential sum
estimates in \cite{hk1} led to $\sigma = 1/7$.  We now describe how to sharpen
the method.

In \cite{hk1}, our estimates for bilinear exponential sums were based partially on a general result
of the second author: Lemma 5.6 in \cite{Kumchev}. However, that result is no longer sufficient
when $\sigma \ge 1/7$. Our first lemma is a variant of Lemma 5.6 in \cite{Kumchev} that can be applied in the present context.
The reader will recognize the major arc arguments from the previous section and note that the complications from
primes in $\mathcal S_q$ continue to be a nuisance.

\begin{lemma}\label{newlemma5.6}
  Suppose that $\alpha$ is real and that $a, q$ are integers with
  \[
    |q\alpha - a| \le Q^2P^{-2}, \quad 1 \le a \le q \le Q^2, \quad (a, q) = 1,
  \]
  and let $\xi_r$ and $\eta_s$ be complex numbers, with $|\xi_r| \le 1$ and $|\eta_s| \le 1$.
  Suppose also that $R, S$ are reals, which together with $z = z(r, s)$  satisfy
  \[
    1 \le R \le V,  \quad 1 \le S \le W, \quad Z \le z(r, s) \le P^{8/35}.
  \]
  Then the exponential sum
  \[
    g(\alpha) = \sum_{r \sim R} \sum_{s \sim S} \sum_{x \sim P/(rs)} \xi_r\eta_s\psi(rsx, z)e(\alpha (rsx)^2)
  \]
  satisfies the inequality
  \[
    g(\alpha) \ll P^{1 + \epsilon}(q + P^2|q\alpha - a|)^{-1/2} + P^{1 - \sigma + \epsilon}.
  \]
\end{lemma}

\begin{proof}
We will initially treat $g(\alpha)$  like the sum $g_j(\alpha)$ studied in the previous section.  We write
\beq\label{3.a}
g(\alpha) = h_1(\alpha) + \sum_{p \in \mathcal S_q} h_p(\alpha) + h^*(\alpha)
\endq
where
$
h_{l}(\alpha)
$
denotes the subsum of $g(\alpha)$ with
$(rsx, q) = l$, and $h^*(\alpha)$ denotes the subsum of $g(\alpha)$
where $(rsx, q) \ge P^{\sigma}$. We note that $h_d(\alpha)$ is a subsum of
$h^*(\alpha)$ unless $d$ is a prime $p$ in the range $Z \le p < P^{\sigma}$.
Of course, if terms involving $p$ arise we must have $q \ge p \ge Z$.  The reader will note that
for $q < z(r,s)$ some of the more awkward terms in the following do not occur.
  We estimate $h^*(\alpha)$ trivially:
  \beq\label{3.b}
    h^*(\alpha) \ll \sum_{ \substack{ d \mid q\\ d \ge P^{\sigma}}} \sum_{ \substack{ k \sim P\\
d \mid k}} k^{\epsilon/2} \ll P^{1 + \epsilon/2} \sum_{ \substack{ d \mid q\\ d \ge P^{\sigma}}} d^{-1} \ll P^{1 - \sigma + \epsilon}.
  \endq

  Now, let $l = 1$ or $l = p \in \mathcal S_q$. We write $q = lq_0$ and $\beta = \alpha - a/q$. We can then use
the arguments of the previous section to write $h_l(\alpha)$ in the form
\[
 \frac 1{\phi(q_0)} \sum_{\chi \bmod q_0} S(\chi, al) \mathop{ \sum_{r \sim R}
\sum_{s \sim S} \sum_{x \sim P/(rs)}}_{l \mid rsx} \xi_r\eta_s\psi(rsx, z)\chi(rsx/l)e(\beta (rsx)^2).
\]
We can then estimate this as the sum of at most three terms of the form
\[
 \frac 1{\phi(q_0)}  \sum_{\chi \bmod q_0} |S(\chi, al)| h_l(\beta, \chi)
\]
with
\[
 h_l(\beta, \chi)  =  \left|W_l^*(\chi, \beta l^2)\right|,
\]
the sum defined by \eqref{GHadd2} with $R' \le R$, $S' \le S$ as before (and so
these parameters satisfy the same inequalities as given in Lemma 3).
When $l = 1$, we deduce that
\beq\label{3.d}
    h_1(\alpha) \ll q^{-1/2 + \epsilon} \sum_{\chi \bmod q}  h_1(\beta, \chi).
  \endq
For $l=p$ we must estimate sums of the form
\beq\label{3.e}
  \frac 1{\phi(q_0)} \sum_{\chi \bmod q_0} |S(\chi, ap)| h_p(\beta, \chi).
  \endq

  Let $q = p^eq_1$, where $e \ge 1$ and $(p, q_1) = 1$. We also write
  \[
    q_2 = \max( qp^{-e}, qp^{-2} ) = \begin{cases} q_1 & \text{if } e = 1, \\ qp^{-2} & \text{if } e \ge 2. \end{cases}
  \]
  We consider two cases.

  Case 1: $e = 1$. Then $(ap, q_0) = 1$ and $q_0 = q_2$, and we have (similarly to the case $l = 1$)
  \beq\label{3.f}
    h_p(\alpha) \ll q_2^{-1/2 + \epsilon} \sum_{\chi \bmod q_2}  h_p(\beta, \chi).
  \endq

  Case 2: $e \ge 2$. Then the exponential sum $S(\chi, ap)$ vanishes when the conductor of $\chi$ is
divisible by $p^{e - 1}$. Otherwise, we have $|S(\chi, ap)| \ll pq_2^{1/2 + \epsilon}$. Hence,
  \beq\label{3.g}
    h_p(\alpha) \ll q_2^{-1/2+\epsilon} \sum_{\chi \bmod q_2} h_p(\beta, \chi\chi_0),
  \endq
  where $\chi_0$ is the principal character modulo $p$. When $e \ge 3$, we have $\chi\chi_0 = \chi$,
so \eqref{3.g} turns into \eqref{3.f}.  When $e = 2$, we deduce from \eqref{3.g} that
  \[
    h_p(\alpha) \ll q_2^{-1/2+\epsilon} \sum_{\chi \bmod q_2} h_p(\beta, \chi) + \Delta
  \]
  where
  \begin{align*}
    \Delta = q_2^{-1/2 + \epsilon}P^{\epsilon/4} \sum_{\chi \bmod q_2} \Big| \sum_{n \sim P/p^2} \theta_n \chi(n) \Big|,
  \end{align*}
  with coefficients $\theta_n$ satisfying $|\theta_n| \le 1$. By Cauchy's inequality and the orthogonality of the characters modulo $q_2$, we obtain
  \begin{align*}
    \Delta^2 &\ll P^{2\epsilon} \sum_{\chi \bmod q_2} \Big| \sum_{n \sim P/p^2} \theta_n \chi(n)
\Big|^2 \ll P^{2\epsilon} \sum_{ \substack{ m, n \sim Pp^{-2}\\ m \equiv n \pmodulo {q_2}}} \phi(q_2) \\
    &\ll P^{2 + 2\epsilon}p^{-4}\big( 1 + pq_2P^{-1} \big) \ll P^{2 + 2\epsilon}p^{-4} \ll P^{2 - 2\sigma + \epsilon}.
  \end{align*}
  Hence,
  \beq\label{3.h}
    h_p(\alpha) \ll q_2^{-1/2 + \epsilon} \sum_{\chi \bmod q_2}  h_p(\beta, \chi)  + P^{1 - \sigma + \epsilon}.
  \endq

  Combining \eqref{3.a}, \eqref{3.b}, \eqref{3.d}, \eqref{3.f} and \eqref{3.h}, we get
  \beq\label{3.i}
    g(\alpha) \ll q^{-1/2 + \epsilon} \sum_{\chi \bmod q}  h_1(\beta, \chi)  + q_2^{-1/2 +
\epsilon}\sum_{\chi \bmod q_2}  h_p(\beta, \chi)  + P^{1 - \sigma + \epsilon}.
  \endq
  By the argument of Lemma 5.1 in \cite{Kumchev},
  \beq\label{3.j}
    \sum_{\chi \bmod q}  h_1(\beta, \chi)  \ll \sideset{}{'} \sum_{\chi \bmod q}  h_1(\beta, \chi)  + h_1(\beta, \chi^0) + \twolinesum{rsx \sim P}{(rsx, q) > 1} \psi(rsx, z),
  \endq
  where $\chi^0$ is the trivial character and $\sum'$ denotes summation over the non-principal characters modulo $q$. Since the last sum on the right side of \eqref{3.j} vanishes when $q< Z$, we have
  \beq\label{3.k}
    \twolinesum{rsx \sim P}{(rsx, q) > 1} \psi(rsx, z) \ll \sum_{\substack{d \mid q\\ d \ge Z}} \sum_{\substack{k \sim P\\ d \mid k}}
k^{\epsilon/4} \ll P^{1 + \epsilon/3}Z^{-1} \ll q^{1/2 - \epsilon}P^{1 + \epsilon}Z^{-3/2}.
  \endq
  We now note that Lemma 5.4 in \cite{Kumchev} remains true if one replaces the hypothesis $z \le P^{23/140}$ of that lemma by the hypotheses
  \[
    z \le P^{8/35}, \quad z\min(R, S) \le P^{11/20}.
  \]
  Thus, the first two terms on the right side of \eqref{3.j} are bounded above by
  \begin{align*}
    P^{1 + \epsilon/2}(1 + P^2|\beta|)^{-1/2} + qP^{11/20 + \epsilon/2}(1 + P^2|\beta|)^{1/2},
  \end{align*}
  provided that
  \[
    z(r, s) \le \min( P^{8/35}, P^{11/20}V^{-1} ) = P^{8/35}.
  \]
  We conclude that
  \beq\label{3.l}
    q^{-1/2 + \epsilon} \sum_{\chi \bmod q} | h_1(\beta, \chi) | \ll
\frac {P^{1 + \epsilon}}{(q + qP^2|\beta|)^{1/2}} + P^{11/20 + \epsilon}Q + P^{1 + \epsilon}Z^{-3/2}.
  \endq

Before estimating the terms arising from $h_p(\alpha)$ we note that there is no contribution unless $z(r,s)\le p \le P^{\sigma}$.
Hence, working in an analogous fashion to the above,
we obtain
  \begin{align}\label{3.m}
    q_2^{-1/2 + \epsilon} \sum_{\chi \bmod q_2}  h_p(\beta, \chi)
&\ll \frac {(P/p)^{1 + \epsilon}}{(q_2 + q_2P^2|\beta|)^{1/2}} + (P/p)^{11/20 + \epsilon}Q + P^{1 + \epsilon}Z^{-2} \notag\\
    &\ll \frac {P^{1 + \epsilon}}{(q + qP^2|\beta|)^{1/2}} + P^{11/20}Q + P^{1 + \epsilon}Z^{-2},
  \end{align}
  provided that
  \[
    W \le (P/p)^{11/20}, \quad z(r, s) \le \min \big( (P/p)^{8/35}, (P/p)^{11/20}V^{-1} \big).
  \]
  These follow from the inequalities
  \[
    1 - 4\sigma \le \frac {11}{20}(1 - \sigma), \quad 3\sigma \le  \frac {11}{20}(1 - \sigma).
  \]
  The desired estimate follows from \eqref{3.i}, \eqref{3.l} and \eqref{3.m}.
\end{proof}

\begin{lemma}\label{HKlemma}
Suppose that $\alpha \in \fm$ and that the function $\rho_j$ in
\eqref{expsums} satisfies hypotheses (iii) and (iv) together with:
\begin{itemize}
  \item [(v)] $\rho_j(m)$ is the linear combination of $O(\cl^c)$
    bilinear sums of the form
    \beq\label{rhodecomp1}
      \sum_{rs = m} \xi_r \eta_s,
    \endq
    where $|\xi_r| \le \tau(r)^c$, $|\eta_s| \le \tau(s)^c$, and either
    $V \le r \le W$, or $r \ge P^{3 \sigma}$ and $\eta_r = 1$ for all $r$.
  \end{itemize}
Then
\beq\label{HKsum}
  f_j(\alpha) \ll P^{1 - \sigma + 2\epsilon}.
\endq
\end{lemma}

\begin{proof}
By Dirichlet's theorem in Diophantine approximation, we can find
integers $a, q$ with
\[
  1 \le q \le (P/Q)^2, \quad (a, q) = 1, \quad
  |q \alpha - a| < (Q/P)^{2}.
\]
Under the assumption of hypothesis (v), the arguments in
Sections 8 and 9 of \cite{Harnew} (see (34) in particular) yield
the bound
\[
  f_j(\alpha) \ll P^{1 - \sigma + \epsilon} +
  P^{1 + \epsilon}\left(\frac{1}{q} + \frac{q}{P^2}\right)^{1/4}.
\]
This establishes \eqref{HKsum} when $q \ge Q^2$. On the
other hand, when $q \le Q^2$ hypotheses (iii) and (iv) ensure
that we can appeal to Lemma \ref{newlemma5.6}. This yields the
bound
\[
  f_j(\alpha) \ll
  \frac {P^{1 + \epsilon}}{(q + P^2|q\alpha - a|)^{1/2}}
  + P^{1 - \sigma + \epsilon},
\]
from which \eqref{HKsum} follows on noting that for $\alpha \in \fm$
we have
\[
  q + P^2|q \alpha - a| > Q.
\]
\end{proof}

\section{The sieve method}
\label{sec5}

We now show how functions $\rho_j$ having properties (i)--(v)
(when $j=2,3$) or (i)--(iv) ($j=4$)
above can be constructed using the sieve method originating
in \cite{Har83} and developed in \cite{Har96,BHP97} by modifying
 the construction used in \cite{Harnew}. Verification of hypotheses (iii) and (iv)
is straightforward, so we shall concentrate on checking hypothesis (v).
It is immediate that $\psi(m,Z)$ satisfies hypothesis
(v) by Theorem 3.1 in \cite{PDS}. Indeed we can
actually obtain the same result for
\beq\label{FT1}
\sum_{r \le V} c_r \psi(m/r,Z),
\endq
where $p|r \Rightarrow p\ge Z$ if $c_r \ne 0$, and $|c_r| \ll 1$.
We now state as a lemma a further refinement.

\begin{lemma}\label{newlem}
Suppose that $p|r \Rightarrow p\ge Z$ if either $c_r \ne 0$
or $b_r \ne 0$ and $|c_r|, |b_r| \le 1$.  Then
\beq\label{newFT}
\twolinesum{r \le V}{s \le Y}
 c_r b_s \psi(m/(rs),Z)
\endq
satisfies hypothesis (v).
\end{lemma}

\begin{proof}
We can reduce the case $rs \le V$ to \eqref{FT1}.
The case $V \le rs \le W$ is immediately
in the correct form. We may therefore suppose that $rs > W$.
Let
\[
\Pi = \prod_{p < Z} p.
\]
Then
\[
\sum_{r,s}  c_r b_s \psi(m/(rs),Z)
=
\sum_{d|\Pi} \mu(d) \sum_{rsnd=m} c_r b_s.
\]
We can then use the technique used in the proof of \cite[Theorem 3.1]{PDS} to
decompose the sum into $O(L^2)$ sums of the requisite types.  The basic idea is to take out the
prime factors of $d$ one by one until a suitable combination of factors lies in the range
from $V$ to $W$ or the size of the ``free variable" exceeds $P^{3 \sigma}$.
This is possible since each of the prime factors is bounded above by $Z$, so
(also using $s \le Y$)
\[
rp_1 \ldots p_u < V \Rightarrow rp_1 \ldots p_{u+1} < P^{1- 4 \sigma} \ \ \text{and} \ \ rsp_1 \ldots p_u < P^{1-3 \sigma}.
\]
This completes the proof.
\end{proof}

Now $\psi(m, P^{1/2})$ is
the characteristic function of the set of primes in $\ci$.
So Buchstab's identity gives
\beq\label{5.1}
\rho_1(m) =  \psi(m, P^{1/2}) = \psi(m, Z) -
  \sum_{Z \le p < P^{1/2}} \psi(m/p,p).
\endq
We first construct $\rho_2(n)$, returning later to \eqref{5.1}
for the lower bound.
Clearly
\[
\begin{split}
\rho_1(n) &\le  \psi(m, Z) -
  \sum_{Z \le p < Y} \psi(m/p,p)
-
\sum_{V \le p < W} \psi(m/p,p)\\
&=\psi_1 - \psi_2 - \psi_3 \ \ \text{say.}
\end{split}
\]
Now $\psi_1$ and $\psi_3$ satisfy hypothesis (v), and
we apply Buchstab's identity again to $\psi_2$:
\[
\psi_2 = \sum_{Z \le p < Y} \psi(m/p,Z) -
\sum_{Z \le q <p < Y} \psi(m/(pq),q)= \psi_4 - \psi_5
\]
say.  Again $\psi_4$ is in the required form and we can apply Buchstab one more time to
$\psi_5$ to obtain:
\[
\psi_5 = \sum_{Z \le q <p < Y} \psi(m/(pq),Z) -
\sum_{Z \le r< q <p < Y} \psi(m/(pqr),r).
\]
The first term on the right hand side above is of the correct form, whereas the second
term can be split into two parts: one which satisfies hypothesis (v), and the rest will be discarded
since it is counted with a negative weight and we are seeking an upper bound.
It follows that the value we obtain for $C_2$, obtained by adding on to $1$ various integrals
corresponding to the discarded sums (compare \cite[Chapter 6]{PDS}), is
\[
\begin{split}
1\ +\ &\int_{1-5 \sigma}^{2 \sigma} w\left(\frac{1-\alpha}{\alpha}\right) \, \frac{d\alpha}{\alpha^2}
\ +\ \int_{1-4\sigma}^{\frac12} w\left(\frac{1-\alpha}{\alpha}\right) \, \frac{d\alpha}{\alpha^2}   \\
&+\ \int_A w\left(\frac{1-\alpha-\beta-\gamma}{\gamma} \right) \, \frac{d\gamma}{\gamma^2}\, \frac{d\beta}{\beta} \, \frac{d\alpha}{\alpha}.
\end{split}
\]
Here $w(u)$ is
Buchstab's function, defined as the continuous solution of
\[
  \begin{cases}
    (uw(u))' = w(u - 1) & \text{if } u > 2, \\
    w(u) = u^{-1}       & \text{if } 1 < u \le 2.
  \end{cases}
\]
Also, $A$ is the three dimensional region given by:
\[
1-6\sigma < \alpha < 1-5 \sigma, \quad 1-6\sigma < \gamma <\beta < \alpha,
\]
with the additional constraint that neither the sum of any pair of variables
nor the sum of all three variables lies in the interval $[2 \sigma, 1-4\sigma]$.
Some simple calculations then yield $C_2 < 1.74$ when $\sigma = 3/20$.

We now begin the construction of our lower bound sieve function by breaking
the final sum on the right hand side of \eqref{5.1} into three parts:
\[
\Psi_1 = \sum_{Z \le p< V} \psi(m/p,p), \quad \Psi_2 = \sum_{V \le p \le W} \psi(m/p,p), \quad
\Psi_3 = \sum_{W < p< P^{\frac12}} \psi(m/p,p).
\]
Hypothesis (v) is met for $\Psi_2$.  Since we cannot cast $\Psi_3$ into a form which satisfies
(v) this term will contribute to $\rho_4(m)$. This contributes $\log(3/2) < 0.406$ to $C_4$.

We apply Buchstab's identity again to $\Psi_1$ thus obtaining
\[
\Psi_1 = \sum_{Z \le p< V} \psi(m/p,Z) - \sum_{Z \le q < p< V} \psi(m/(pq),q) = \Psi_4 - \Psi_5
\]
say.
By \eqref{FT1} $\Psi_4$ satisfies (v).  We split $\Psi_5$ into four sums $\Sigma_j, 1 \le j \le 4$ according to
the sizes of $p, q$ as follows:
\[
\begin{split}
&j=1: \qquad V \le pq \le W;\\
&j=2: \qquad pq>W, q> Y;\\
&j=3: \qquad pq> W, q \le Y;\\
&j=4: \qquad pq<V.
\end{split}
\]
Now $\Sigma_1$ automatically satisfies (v).  We must discard the whole of $\Sigma_2$ and this is the main contribution to
$\rho_5(n)$ leading to a ``loss" at $\sigma=3/20$:
\[
\int_{1/4}^{3/10} \int_{1/4}^{\alpha} \frac{d\beta \,  d \alpha }{\alpha \beta (1 - \alpha - \beta)} < 0.037.
\]
where we have noted that $w(u) = 1/u$ throughout the integration region. We can apply Buchstab's identity again to $\Sigma_3$, leading to
\[
\sum_{pq>W, p < V}\psi(m/(pq),Z)     - \sum_{pq>W, Z<r<q<Y} \psi(m/(pqr),r).
\]
The first term here can be treated using Lemma \ref{newlem}.  The second term can be split into three more sums depending on
whether: $V \le pr \le W$, in which case (v) is automatically satisfied; $pr > W, qr>Y$ in which case we discard this portion
which leads to another term in $\rho_4$ whose contribution to $C_4$ at $\sigma=3/20$ is $< 0.08$; $qr < Y$ in which case we can decompose once more and the resulting sums all satisfy hypothesis (v) since
\[
Z\le s<r<q, rq<Y \Rightarrow V \le qrs \le W
\]
when $\sigma \le 3/20$.

Finally, we can apply Buchstab's identity again to $\Sigma_4$ to obtain
\[
\twolinesum{Z \le q < p< V}{pq < V} \psi(m/(pq),Z) - \twolinesum{Z \le r < q < p< V}{pq < V} \psi(m/(pqr),r).
\]
The first sum above satisfies hypothesis (v) by Lemma \ref{newlem}.  We split the second sum
into two sums, one with $V \le pqr \le W$ and one with $pqr>W$.  The first sum
immediately satisfies (v), while we can apply Buchstab's identity to the second since $pq < V, r <q<V^{1/2}<Y$.
This leads to a sum over four prime variables which can often be grouped into products lying between
$V$ and $W$.  The rest of this sum leads to one last contribution to $\rho_5$ from a four dimensional integral whose contribution
at $\sigma = 3/20$ is $<0.0006$.

We can now gather all our results together to obtain
\[
C_3 - C_2 C_4 > (1+C_4 - 0.038) - 1.74 C_4 > 0.96 - 0.74 \times 0.49 = 0.5974,
\]
as needed to complete the proof.

\section{Proof of Theorem 1}
The proof follows a very similar pattern to our previous paper, but we include all the details for completeness.
There is one additional complication we must first deal with that did not arise in our earlier work.  That is,
we need a bound for
\[
\int_0^1 |g_j(\alpha)|^4 \, d \alpha \quad \text{and not just for} \ \ \int_0^1 |f_j(\alpha)|^4 \, d \alpha.
\]
The bound $\ll P^{2 + \epsilon}$ follows for the latter integral immediately from Hua's lemma (Lemma~2.5 in \cite{Vaughan}),
but the result we require demands a little more work.

\begin{lemma}\label{newhua}
In the notation of previous sections, for $j=2$ or $4$ we have, for any $\epsilon > 0$,
\[
\int_0^1 |g_j(\alpha)|^4 \, d \alpha \ll P^{2 + \epsilon}.
\]
\end{lemma}

\begin{proof}
We recall the set $\mathcal D = \{m: P^{\sigma} \le m \le Q, p|m \Rightarrow p > Z \}$.  Further, put
$\mathcal D_q = \{d \in \mathcal D : (d,q) \ge P^{\sigma}\}$. We then have
\[
\begin{split}
\int_0^1 |g_j(\alpha)|^4 \, d \alpha
&\ll
\int_0^1 |f_j(\alpha)|^4  \, d \alpha+  \int_0^1 |g_j(\alpha)- f_j(\alpha)|^4 \, d \alpha\\
&
\ll P^{2 + \epsilon} + \sum_{P^{\sigma} \le q \le Q} \sum_{(a,q) = 1} I(a,q),
\end{split}
\]
where
\[
I(a,q) = \int_{\fM(q,a)} \Bigg|\sum_{d \in {\mathcal D}_q} \twolinesum{m \in \ci}{d|m}\rho_j(m) e(\alpha m^2)     \Bigg|^4 \, d\alpha,
\]
since it is only on these arcs that $\theta(m,\alpha) = 0$. By H\"older's inequality
\[
I(a,q) \le \left(\sum_{d \in {\mathcal D}_q} 1 \right)^3  \sum_{d \in {\mathcal D}_q} \int_{\fM(q,a)} \Bigg|\twolinesum{m \in \ci}{d|m}\rho_j(m) e(\alpha m^2)\Bigg|^4
\, d\alpha.
\]
Since each $d \in \mathcal D_q$ has at most two prime divisors a simple change of integration variable
and a swap in the order of summation gives
\[
\sum_{P^{\sigma} \le q \le Q} \sum_{(a,q) = 1} I(a,q)
\ll
\sum_{d \in \mathcal D} \frac{1}{d^2} \twolinesum{q \le Q}{(d,q) > P^{\sigma}}\sum_{(a,q) = 1}
\int_{\fM^*(a,d,q)} |\Sigma(\alpha,d)|^4 \, d\alpha,
\]
where
\[
\fM^*(a,d,q) = \left[\frac{ad^2}{q} - \frac{\omega d^2}{q}, \frac{ad^2}{q} + \frac{\omega d^2}{q}\right]
\]
and
\[
\Sigma(\alpha,d) = \sum_{md \in \ci} \rho_j(md) e(\alpha m^2).
\]
Since $\omega d^2/q < Q^{-2}$ the intervals $\fM^*(a_1,d,q), \fM^*(a_2,d,q)$ overlap $\pmodulo 1$ only if
$a_1d^2 \equiv a_2 d^2 \pmodulo q$. The number of overlaps is thus $\le (q,d^2)$
as $a$ runs through the reduced residues $\pmodulo q$.  Now, since $(q,d) \ge P^{\sigma}$
and $q< P^{2 \sigma}$, we must have $(q,d^2) = (q,d)$.  This value is either $d$ or otherwise the larger of the
two prime divisors of $d$, which we denote by $p(d)$ if $d$ is not a prime (and let $p(d)$ be zero if $d$ is prime).
Now, by Hua's lemma
\[
\int_{\fM^*(a,d,q)} |\Sigma(\alpha,d)|^4 \, d\alpha \le \int_0^1 |\Sigma(\alpha,d)|^4 \, d\alpha \ll \left(\frac{P}{d}\right)^{2 + \epsilon}.
\]
Hence
\[
\begin{split}
\sum_{P^{\sigma} \le q \le Q} \sum_{(a,q) = 1} I(a,q)
&\ll
P^{2 + \epsilon} \sum_{d \in \mathcal D} \frac{1}{d^{4+ \epsilon}} \left(\sum_{q' \le Q/d} d + \sum_{q' \le Q/p(d)} p(d)\right) \\
&\ll P^{2 + \epsilon} Q \sum_{d \in \mathcal D} \frac{1}{d^{4}} \ll \frac{P^{2 + \epsilon} Q}{P^{3 \sigma}} < P^2.
\end{split}
\]
This completes the proof.
\end{proof}

Let $\fz$ be the set of integers $n \in \cb$ for which
\eqref{MA} holds but which are not representable as sums of
three squares of primes. We write $|\fz|$ for the cardinality
of $\fz$ and $Z(\alpha)$ for its generating function:
\[
  Z(\alpha) = \sum_{n \in \fz} e(-\alpha n).
\]
Write
\[
\begin{split}
G^*(\alpha)
&= f_1(\alpha)\left(f_1(\alpha) f_3(\alpha) - f_2(\alpha)f_4(\alpha)\right),\\
G(\alpha)
&= f_1(\alpha)\left(f_1(\alpha) g_3(\alpha) - g_2(\alpha)g_4(\alpha)\right),\\
K(\alpha)
&=
G^*(\alpha) - G(\alpha)\\
&=
f_1^2(f_3 - g_3) + f_1\left(f_2(f_4-g_4) + g_4(f_2 - g_2)\right),
\end{split}
\]
where we have omitted the common variable $\alpha$ for all the
functions on the last line in the interests of clarity. Then
\[
  \int_0^1 G^*(\alpha) Z(\alpha) \, d\alpha \le 0
\]
and
\[
  \int_{\fM} G(\alpha) Z(\alpha) \, d \alpha
  \gg |\fz|P\cl^{-6}.
\]
Thus,
\[
  \bigg| \int_{\fm} G(\alpha) Z(\alpha) \, d\alpha +
  \int_0^1 K(\alpha) Z(\alpha) \, d\alpha \bigg|
  \gg |\fz|P\cl^{-6}.
\]
Recalling Lemma \ref{HKlemma} and \eqref{2.5}, we deduce that
\begin{align*}
  |\fz| &\ll L^6P^{-1} \bigg(
  \int_{\fm} \big| G(\alpha) Z(\alpha) \big| \, d\alpha
  + \int_0^1 \big| K(\alpha)Z(\alpha) \big| \, d\alpha
  \bigg) \\
  &\ll P^{-3/20+\epsilon/2}
  \int_0^1 \big| h(\alpha)^2 Z(\alpha) \big| \, d\alpha,
\end{align*}
where $h(\alpha)$ is one of the $f_j(\alpha), g_j(\alpha)$.
Finally, using Cauchy's inequality, Parseval's identity and Hua's lemma (or Lemma \ref{newhua} if $h(\alpha) = g_j(\alpha)$ for some $j$),
we find that the last integral is bounded by
\[
  \left(\int_0^1 |Z(\alpha)|^2 d\alpha \right)^{1/2}
  \left(\int_0^1 |h(\alpha)|^4 d\alpha \right)^{1/2}
  \ll |\fz|^{1/2}P^{1 + \epsilon/2},
\]
and so
\[
  |\fz| \ll P^{17/10 + 2\epsilon} \ll N^{17/20 + \epsilon}.
\]
 Combining this estimate with Lemma \ref{MAlem} then proves Theorem 1 as required. \qed

\bibliographystyle{amsplain}

\end{document}